\documentclass{article}
\usepackage{amsmath}
\usepackage{graphicx}
\usepackage{verbatim} 
\newcommand{\qed}{{\unskip\nobreak\hfil\penalty50\hskip2em\vadjust{}
            \nobreak\hfil$\Box$\parfillskip=0pt\finalhyphendemerits=0\par}}

\newtheorem{lemma}{Lemma}[section] 
\newtheorem{corollary}{Corollary}[section]

\newtheorem{definition}{Definition}[section]
\newcommand{\bed}{\begin{definition}}
\newcommand{\eed}{\end{definition}}

\newtheorem{theorem}{Theorem}

\newcommand{\eps}{\epsilon}

\newcommand{\bitem}{\begin{itemize}}
\newcommand{\eitem}{\end{itemize}}

\newcommand{\goto}{\rightarrow}

\newcommand{\beqn}{\begin{equation}}
\newcommand{\eeqn}{\end{equation}}
\newcommand{\balign}{\begin{align}}
\newcommand{\ealign}{\end{align}}

\newcommand{\hL}{\hat{L}}

\newcommand{\pee}{\pi}
\newcommand{\Pee}{P}
\newcommand{\wts}{w}
\newcommand{\hCov}{\widehat{Cov}_{n,p}}
\newcommand{\hwts}{\hat{w}}

\usepackage{amssymb}
\begin{document}

\title{Feature Selection by  Higher Criticism Thresholding:\\
Optimal Phase Diagram}
\author{David Donoho$^1$ and Jiashun Jin$^2$ \\
$^1$Department of Statistics, Stanford University\\
$^2$Department of Statistics,     Carnegie Mellon  University}
\date{December  11, 2008}
\maketitle

\begin{abstract}
We consider two-class 
linear classification in a high-dimensional, low-sample size setting. 
Only a small fraction of the features are useful, 
the useful features are unknown to us, and each useful feature 
contributes weakly to the classification decision --
this setting was called the rare/weak model (RW Model) in \cite{PNAS}.

We select features by thresholding feature $z$-scores.  
The  threshold is set by {\it higher criticism} (HC)  \cite{PNAS}.
Let $\pee_i$ denote the $P$-value associated to the $i$-th $z$-score and
$\pee_{(i)}$ denote the $i$-th order statistic of the collection of $P$-values.
The HC threshold (HCT) is the order statistic of the $z$-score corresponding to index $i$
maximizing  $(i/n - \pee_{(i)})/\sqrt{\pee_{(i)}(1-\pee_{(i)})}$.  
The ideal threshold optimizes the classification error. 
In \cite{PNAS} we showed that HCT
was numerically close to the ideal threshold.

We formalize an asymptotic
framework for studying the RW model, considering
a sequence of problems with increasingly many features and relatively fewer observations. 
We show that along this sequence,
the limiting performance of ideal HCT is essentially just as good as 
the limiting performance of ideal thresholding.
Our results describe two-dimensional {\it phase space},
a two-dimensional diagram
with coordinates  quantifying ``rare" and ``weak" in the RW model.
Phase space can be partitioned into two regions --
one where ideal threshold classification is successful, and one
where the features are so weak and so rare that it must fail. 
Surprisingly, the regions where ideal HCT 
succeeds and fails make the exact same
partition of the phase diagram. 
Other threshold methods,
such as FDR threshold selection,
are successful in a substantially smaller region
of the phase space than either HCT or Ideal thresholding.

The False Feature Detection Rate (FDR) and local  FDR 
of the ideal and HC threshold selectors have surprising
phase diagrams which we also describe.
Results showing asymptotic equivalence
of HCT with  ideal HCT 
can be found 
in a companion paper \cite{AoS}.
\end{abstract}
\begin{quote} \small
{\bf Keywords}:  Asymptotic rare weak model (ARW),   False Discovery Rate (FDR),    linear classification,  local False Discovery Rate (Lfdr),   phase diagram,  Fisher's separation measure (SEP),  feature selection by thresholding.
\end{quote}
\begin{quote} \small
{\bf AMS 2000 subject classifications}:  Primary-62H30, 62H15.
Secondary-62G20, 62C32, 62G10.
\end{quote}

{\bf Acknowledgments}:   The authors would like to thank
the Isaac Newton Mathematical Institute of Cambridge University
 for hospitality during the program
{\it Statistical Theory and Methods for Complex, High-Dimensional Data},  
and for a Rothschild Visiting Professorship held by DLD.  DLD was partially 
supported by NSF DMS 05-05303,  and  JJ was partially supported by NSF CAREER award  DMS 06-39980.    
\newcommand{\ntrain}{n}
\section{Introduction}
\setcounter{equation}{0} 

The modern era of high-throughput data collection
creates data in abundance.
Some devices -- spectrometers and gene chips come to mind --
automatically  generate measurements
on thousands of standard features of each observational unit.

What hasn't changed in science is the difficulty of obtaining
good observational units.   High-throughput devices don't help
us to find and enroll qualified patients, or catch exotic butterflies,
or observe primate mating behaviors.  Hence, in many fields 
the number of observational units -- eg. patients, butterflies, or matings --
has not increased, and stays today  in the dozens or hundreds. 
But each of those few observational units
can now be subjected to a large battery of
automatic feature measurements

Many of those automatically measured features
will have little relevance to any given project.  
This new era of {\it feature glut} poses
a needle-in-a-haystack problem: we must detect a relatively few valuable
features among many useless ones. 
Unfortunately,  the combination of
small sample sizes (few observational units)
and high dimensions (many feature measurements) makes it hard to tell
needles from straw.

Orthodox statistical methods assumed a quite
different set of conditions:  more observations
than features, and all features highly relevant.
Modern statistical research
is intensively developing new tools
and theory to address the new unorthodox setting; such
research comprised much of the activity  in the 
recent 6-month Newton Institute program 
{\it Statistical Theory and Methods for Complex,   High-Dimensional Data}.

In this article we focus on this new setting, this time
addressing the challenges that modern
high-dimensional data pose to linear classification
schemes. New data analysis tools and
new types of mathematical analysis of those tools
will be introduced.
 
\subsection{Multivariate normal classification}
\label{sec-simple}

Consider a simple model of linear classifier training.
We have a set of labelled training samples $(Y_i, X_i)$, $i=1,\dots,n$,
where each label $Y_i$ is $ \pm 1$ and each feature vector $X_i \in R^p$.
 For simplicity, we  assume the training set contains equal numbers of $1$'s and
$-1$'s and that the feature vectors $X_i \in R^p$ obey
$X_i \sim N( Y_i \mu , \Sigma )$,  $i=1,\dots, n$,
for an unknown {\it mean contrast  vector} $\mu \in R^p$; here
$\Sigma$ denotes the feature covariance matrix and $\ntrain$ is
the training set size. In this simple setting, one ordinarily 
uses  linear classifiers   to classify an unlabeled test vector $X$, taking the general  form $L(X)  =  \sum_{j=1}^p \wts(j) X(j)$, 
for  a sequence of `feature weights' $\wts = (\wts(j): j=1,\dots,p)$.

Classical theory going back to RA Fisher \cite{Anderson} shows that the 
optimal classifier has feature weights $\wts \propto \Sigma^{-1} \mu$;
at first glance linear classifier design seems
straightforward and settled.
However,  in many of today's most active application areas,
it is a major challenge to construct linear classifiers
which work well.  

\subsection{$p$ larger than $n$}
\label{sec-p-gt-n}
The culprit can be called the ``$p > n$ problem''.
A large number $p$ of measurements is automatically made
on thousands of standard features, but in a given project,
the number of observational units, $\ntrain$,
might be in the dozens or hundreds. 
The fact that $p \gg n$ makes it
difficult or impossible to 
estimate the feature covariance $\Sigma$ 
with any precision.

It is well known that naive application of the formula 
$\wts \propto \Sigma^{-1} \mu$ to empirical data in the $p > n$ setting
is problematic; at a minimum, because
the matrix of empirical feature covariances $\hCov(X)$ 
is not invertible. But even if we use the generalized inverse $\hCov(X)^\dag$,
the resulting naive classification weights, 
$\hwts_{naive}\propto  \hCov(X)^{\dag} \hCov(Y,X)$,
often give very ``noisy'' classifiers with low accuracy.
The modern feature glut
 thus seriously  damages the applicability of
 `textbook' approaches.
 
A by-now standard response to this problem
is to simply ignore feature covariances, 
and standardize the features to
mean zero and variance one.   One in effect pretends that
the feature covariance matrix $\Sigma$ is the identity matrix,
and uses the formula
$\wts(j) \propto Cov(Y,X(j)) $ \cite{Bickel,Fan}.
Even after this reduction, further challenges remain.

\subsection{When features are rare and weak}   \label{sec-simple-model} 

In fields such as genomics, with
a glut of feature measurements per observational unit,
it is expected
 that few measured features will be useful in any given project;
 nevertheless, they all get measured, because researchers can't say
in advance which ones will
be useful. Moreover,
reported misclassification rates are relatively high.
Hence there may be numerous useful features,
but they are relatively rare and individually quite weak. 

Such thinking motivated   the following {\it rare/weak feature model} {\it (RW Feature Model)}
in \cite{PNAS}.
Under this model, 
\bitem
\item {\bf Useful features are rare:} the contrast vector $\mu$ is
nonzero in  only $k$ out of $p$ elements, where
 $ \eps = k/p$ is small, i.e. close to zero. As an example,
think of $p = 10,000$,    $k =100$, so $\eps = k/p = .01$.
In addition, 
\item  {\bf Useful features are weak:} the nonzero elements
of $\mu$ have {\it common} amplitude $\mu_0$,
which is assumed not to be `large'. `Large' can be measured 
using $\tau = \sqrt{n} \mu_0$; values of $\tau$
in the range $2$ to $4$ imply that corresponding values of $\mu$ are
not large. 
\eitem 
Since the
elements $X(j)$ of the feature vector where 
the class contrast $\mu(j) = 0$
are entirely uninformative about the value of $Y(j)$, 
only the $k$ features where $\mu(j) = \mu_0$
are useful.  The problem is how to identify and benefit from
those rare, weak features.  

We speak of 
$\eps$ and $\tau$ as the {\it sparsity} and {\it strength} parameters
for the Rare/Weak model, and refer to the $RW(\eps,\tau)$ model.
Models with a `sparsity' parameter $\eps$ are common 
in estimation settings  \cite{DJ,MENBO,Needles}, but not with 
the feature strength constraint $\tau$.
Also closely related to the RW model is
work in multiple testing by Ingster and the authors \cite{Ingster97,HC,JinThesis}.

\subsection{Feature selection by thresholding}

Feature selection \-- 
i.e. working only with an empirically-selected
subset of features \-- is a standard response to feature glut.  
We are supposing, as announced in 
Section \ref{sec-p-gt-n}, that
feature correlations can be ignored and 
that features are already standardized to variance one.
We therefore focus on 
the vector of feature $Z$-scores,
with components $ Z (j) = n^{-1/2} \sum_i Y_i X_i(j)$,  $j=1,\dots,p$.
These are the $Z$-scores of two-sided  tests
of $H_{0,j}$:  $Cov(Y,X(j)) = 0$.
Under our assumptions
  $Z \sim N(\theta , I_p)$ with $\theta=\sqrt{n} \mu$
  and $\mu$ the feature contrast vector.
Features with nonzero $\mu(j)$
typically have significantly nonzero $Z(j)$ while all
other features will have $Z(j)$ values largely consistent
with the null hypothesis $\mu(j)= 0$. 
In such a setting, selecting features with $Z$-scores above
a threshold makes sense. 

We identify three useful threshold functions:
 $\eta_t^{\star}(z)$,   $\star \in \{ clip, hard,soft\}$. These are:
 {\it Clipping} -- $ \eta^{clip}_t(z) = \mbox{sgn}(z) \cdot 1_{\{|z| > t\}}$, which ignores  the size of the $Z$-score, provided it is large;   {\it Hard Thresholding} -- 
 $ \eta^{hard}_t(z) =  z\cdot 1_{\{ |z| > t \}}$, which uses
 the size of the $Z$-score, provided it is large; and 
  {\it Soft Thresholding} -- 
 $ \eta^{soft}_t(z) =   \mbox{sgn}(z) (|z| - t)_+$, which uses
a shrunken $Z$-score, provided it is large.
\begin{definition}
Let  $\star \in \{ soft, hard, clip\}$.
The threshold feature selection classifier makes its decision 
based on  $L_t^\star <> 0$ 
where  $ \hL_t^\star(X) = \sum_{j=1}^p \hwts^\star_t(j) X(j)$,   and 
$ \hwts^{\star}_t(j) = \eta^\star_t ( Z(j)),  j = 1,\dots,p$. 
\end{definition}
In words, the classifier sums across features with
large training-set $Z$-scores, and a simple function of the $Z$-score
generates the feature weight.

Several methods for linear classification in bioinformatics
follow this approach: the Shrunken Centroids method 
\cite{Tibs} is a variant of soft thresholding in this two-class setting; 
the highly-cited methods in  \cite{Leukemia} and \cite{BreastCancer} 
are variants of hard thresholding. Clipping makes sense
in the theoretical setting of the RW model (since then the useful features
have all the same strength)  and is simpler to analyse
than the other nonlinearities.

Thresholding has been popular in estimation
for more than a decade \cite{DJ}; it is known to be succesful
in `sparse' settings where the estimand has many coordinates,
of which only a relatively few coordinates are significantly nonzero.
However classification is not the same as estimation, and performance
characteristics are driven by quite different considerations.

One crucial question remains: how to choose the threshold based on the data? 
Popular methods for threshold choice include cross-validation
\cite{Tibs};    control of the false discovery rate    \cite{AB,ABDJ,BH95,FDR}; 
and control of the local false discovery rate  \cite{Efron}.  

\subsection{Higher Criticism}

In \cite{PNAS} we proposed a method of threshold choice based on
recent work in the field of multiple comparisons. We now very briefly 
mention work in that field and then introduce the threshold choice method.

\subsubsection{HC testing}
Suppose we have a collection of $N$ $\Pee$-values
$\pee_i$ which under the global null hypothesis
are uniformly distributed: $\pee_i \sim_{iid} U[0,1]$.
Consider the order statistics:
$ \pee_{(1)}  \leq \pee_{(2)} \leq \dots \leq \pee_{(N)}$.
Under the null hypothesis,
these order statistics have the usual properties of
uniform order statistics, including the asymptotic normality
$\pee_{(i)} \sim_{approx}  \mathrm{Normal}( i/N, i/N ( 1 - i/N))$.
HC forms a $Z$-score comparing $\pee_{(i)}$ with its
mean under the null, and then maximizes over a wide range of $i$.
Formally:
\begin{definition} ({\bf HC Testing}) \cite{HC}. 
The Higher Criticism objective    is 
\begin{equation} \label{HCobj}
HC(i; \pee_{(i)}) = \sqrt{N} \frac{i/N - \pee_{(i)}}{\sqrt{i/N(1-i/N)}} .
\end{equation}
 Fix $\alpha_0 \in (0,1)$ (eg $\alpha_0 = 1/10$).  
 The HC test statistic is  $HC^* = \max_{1 \leq i \leq \alpha_0 N } HC(i; \pee_{(i)})$.
\end{definition}

HC seems insensitive to the selection of $\alpha$, in Rare/Weak situations;
here we always use $\alpha_0 = .10$.

In words, we look for the largest standardized discrepancy 
for any $\pi_{(i)}$  between the observed behavior and the expected
behavior under the null.   When this is
large, the whole collection of $\Pee$-values is not consistent with 
the global null hypothesis.  The phrase 
``Higher Criticism'' is due to John Tukey, and
reflects the shift in emphasis from
single test results to the whole collection of tests;
see discussion in  \cite{HC}.  
 Note:  there are several variants of HC statistic; we discuss only one variant 
in this brief note; the main results of \cite{HC} still apply to this variant;
for full discussion see \cite{HC,PNAS}.

\subsubsection{HC thresholding}
\label{ssec-HCTDEF}

Return to the classification
setting 
in previous sections. 
We have  a vector of feature $Z$-scores  $(Z(j), j=1,\dots, p)$.
To apply HC notions, translate 
$Z$-scores into two-sided $\Pee$-values, and maximizes  the HC
objective over index $i$ in the appropriate range.  Define the {\it feature $\Pee$-values}
$\pee_i = Prob \{ |N(0,1)| > |Z(i)|\}$,   $i=1,\dots, p$;
and define the  increasing rearrangement $\pee_{(i)}$,
the HC objective function $HC(i; \pee_{(i)})$, and the increasing rearrangement $|Z|_{(i)}$ 
correspondingly.
Here is our proposal.
\begin{definition} ({\bf HC Thresholding}).   \label{hct-def}
Apply the HC procedure to the feature $\Pee$-values.
Let the maximum HC objective be achieved at index $\hat{i}$.
The {\bf Higher Criticism threshold} (HCT) is the value $\hat{t}^{HC} = |Z|_{(\hat{i})}$.
The {\bf HC threshold feature selector} selects features with
$Z$-scores exceeding $\hat{t}^{HC}$ in magnitude.
\end{definition}

\begin{figure} 
\begin{centering}
\includegraphics[height=2.75 in,width=3.15in]{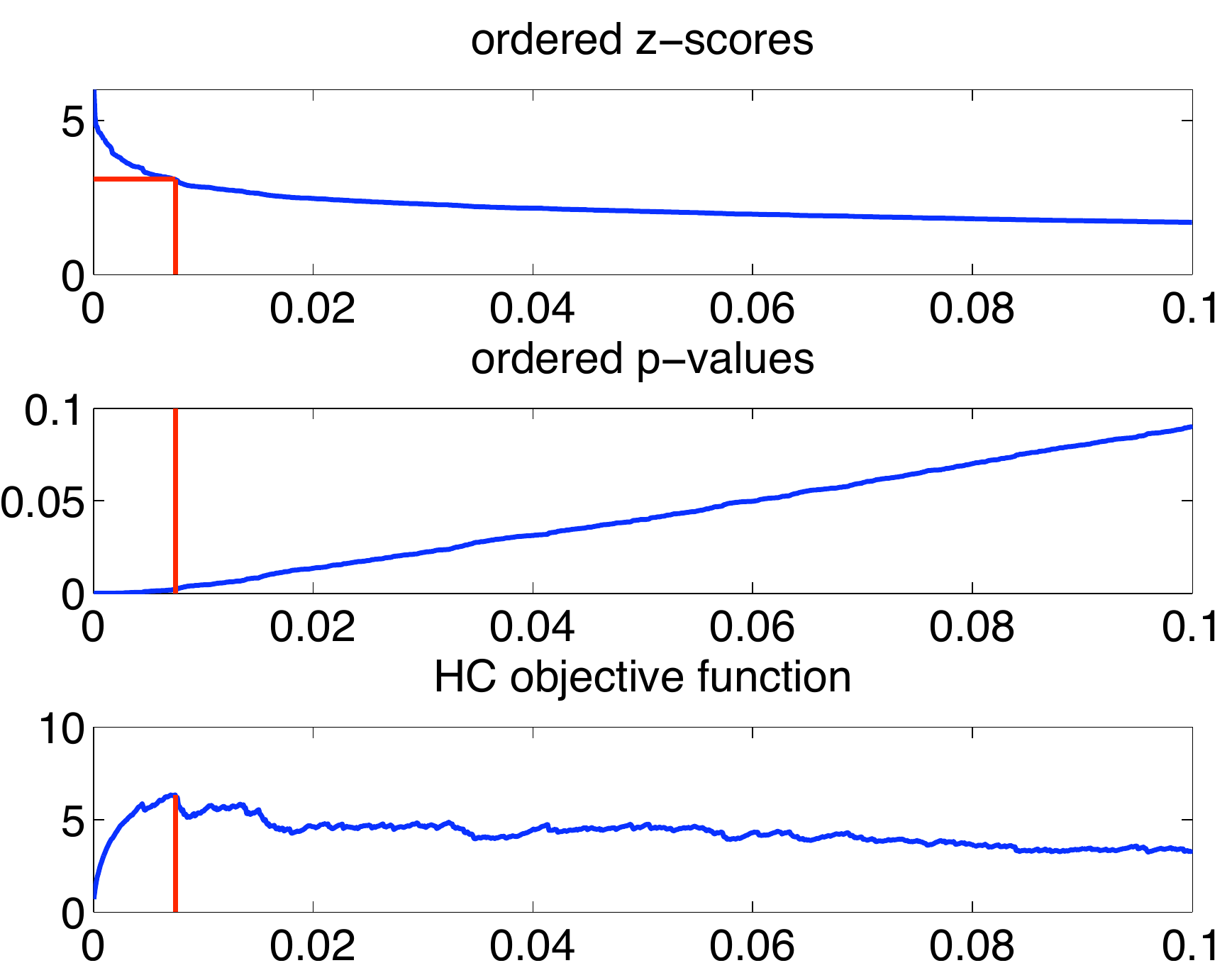}
\caption{Illustration of HC thresholding. Panel (a) the ordered $|Z|$-scores.
Panel (b) the corresponding ordered $\Pee$-values in a PP plot. Panel (c)
the HC objective function in    (\ref{HCobj}); this is largest at
$\hat{i} \approx 0.007 p $ ($x$-axes are $i/p$).   Vertical  lines indicate $\pi_
{(\hat{i})}$ in Panel (b), and $|Z|_{(\hat{i})}$
in Panel (a). }
\label{fig:IllustHC}
\end{centering}
\end{figure}

Figure     \ref{fig:IllustHC}    illustrates the procedure.
Panel (a) shows a sample of $Z$-scores,  Panel (b) shows
a PP-plot of the corresponding ordered $\Pee$-values versus $i/p$ and 
Panel (c) shows a standardized PP-plot.
The standardized PP-Plot has its largest deviation from zero at $\hat{i}$;
and this generates the threshold value.

\subsubsection{Previously-reported results for HCT}

Our article  \cite{PNAS} reported several findings about behavior of HCT
based on numerical and empirical evidence.
In the RW model,  we can define an ideal threshold,
i.e. a threshold based on full knowledge of 
the RW parameters $\eps$ and $\mu$ and chosen to minimize the
misclassification rate of the threshold classifier -- see Section   \ref{sec-ideal} below.
We showed in \cite{PNAS} that:
\bitem
 \item HCT gives a threshold value which is numerically
 very close to the ideal threshold.
  \item In the case of very weak feature z-scores,
 HCT has a False Feature Discovery Rate (FDR) substantially higher
 than other popular approaches, but a Feature Missed Detection Rate (MDR)
 substantially lower than those other approaches;
 \item At the same time, HCT has FDR and MDR very closely
 matching those of the ideal threshold.
 \eitem
In short, HCT has very different operating characteristics
from those other thresholding schemes like FDR
thresholding \cite{AB,ABDJ} 
 and Bonferroni thresholding, but very similar operating characteristics
to the ideal threshold.

\subsection{Asymptotic RW model, and the phase  diagram}

In this paper we further support
the findings reported in \cite{PNAS},
this time using an asymptotic analysis.
In our  analysis
the number of observations $n$ and 
the number of features $p$ tend to infinity in a linked fashion,
with $n$ remaining very small compared to $p$. (Empirical results in \cite{PNAS} show
 our large-$p$ theory is applicable at moderate $n$ and $p$).

More precisely, we consider a sequence of problems with increasingly
more features, increasingly more rare useful features, and relatively small
numbers of observations compared to the number of features.

 \begin{definition}
 \label{def-ARW}
The phrase {\bf asymptotic RW model} 
refers to the following combined assumptions.
\bitem
\item {\bf Asymptotic Setting.} We consider a sequence of problems,
where the number of observations $n$ and the number of features $p$
both tend to $\infty$ along the sequence.
\item 
{\bf $p$ dramatically larger than $n$}.
Along this sequence, 
$n \sim c \cdot \log(p)^\gamma $, so there are dramatically more
features per observational unit than there are observational units. 
\item
{\bf  Increasing  Rarity.} The sparsity  $\eps$ varies
with $n$ and $p$ according to  $\eps = p^{-\beta}$, $0 < \beta < 1$.
\item {\bf Decreasing Strength.} The  strength  $\tau$ varies with $n$ and $p$ 
 according to $\tau = \sqrt{2r \log(p)}$, $0 < r < 1$.
\eitem
The symbol $ARW(r,\beta;c,\gamma)$ refers to the model combining these
 assumptions.
\end{definition}
In this model, because $r < 1$,
useful features are individually too weak to detect
and because $0 < \beta < 1$,
useful features are increasingly rare with increasing $p$,
while increasing in total number with $p$. 
It turns out that $c$ and $\gamma$ are incidental,
while $r$ and $\beta$ are the driving parameters. 
Hence we always simply write $ARW(r,\beta)$ below.

There is a large family of choices of $(r,\beta)$
where successful classification is possible, and another
large family of choices where it is impossible.
To understand this fact,
we use the concept of
 {\it phase space}, the two-dimensional
domain $0 < r , \beta < 1$.
We show  that this domain 
is partitioned into two regions or `phases'.
In the ``impossible'' phase, useful features are so
rare and so weak that  classification
is asymptotically impossible {\it even with the ideal choice
of threshold}.  In the ``possible''
phase,  successfully separating the two groups is
indeed possible \-- {\it if} one has access to 
the ideal threshold.  Figure \ref{fig:phase}
displays this domain and its partition into phases.
Because of the partition into two phases, we also
call this display the {\it phase diagram}.
An explicit formula for the graph $r  = \rho^*(\beta)$ 
bounding these phases is given in (\ref{eq-rho-def}) below.

{ The phase diagram provides a convenient platform for comparing 
different procedures.     A threshold choice is {\bf optimal} 
if it gives the same partition of phase space as the
one obtained with the ideal choice of threshold}.  

How does HCT compare to the ideal threshold,  and what  partition in the phase space   does HCT  yield? 
For reasons of space, we focus in this paper
on the {\it Ideal HC threshold}, which is obtained upon replacing the empirical
distribution of feature $Z$-scores by it expected value. The Ideal
HC threshold is thus the threshold which HCT is `trying to estimate';
in the companion paper \cite{AoS}  we give a full analysis showing that the 
ideal HCT and HCT are close.

The central surprise of our story is that HC behaves surprisingly well:   
the partition  of phase space  
describing the two regions where ideal thresholding
fails and/or succeeds also describes
the two regions where Ideal HCT fails and/or succeeds
in classifying accurately. 
The situation is depicted in the table below:
\begin{table}[h]
\begin{center}
\begin{tabular}{|l|l|l|}
\hline
Region & Property of Ideal Threshold & Property of Ideal HCT \\
\hline
$r < \rho^*(\beta)$ & Ideal Threshold Classifier Fails & Ideal HCT Fails \\
$r > \rho^*(\beta)$ & Ideal Threshold Classifier Succeeds & Ideal HCT Succeeds \\
\hline
\end{tabular}
\end{center}
\end{table}
Here by `succeeds', we mean asymptotically zero misclassification rate
and by  `fails', we mean asymptotically 50\% misclassification rate.
 
In this sense of size of regions of success, HCT
is just as good as the ideal threshold.  
Such statements cannot be made for some other popular thresholding schemes,
such as False Discovery threshold selection. 
As will be shown in \cite{AoS} even the very popular 
Cross-Validated choice of Threshold will fail if the training set size is bounded,
while HCT will still succeed in the RW model in that case.

The full proof of a broader set of claims 
-- with a considerably more general treatment --
will appear elsewhere.     { To elaborate 
the whole story on HCT needs three connected papers including   \cite{PNAS},   \cite{AoS}, and the current one.  
In \cite{PNAS},     we  reported numerical results both with simulated 
data from the RW model  and with certain real data often used 
as standard benchmarks for classifier performance.   
In \cite{AoS},    we will develop a more mathematical treatment of 
many results we cite here and in \cite{PNAS}.   
The current article, logically second in  the triology, 
develops an analysis of Ideal HCT which is both
transparent and which provides the key insights
underlying our lengthy arguments in \cite{AoS}}.     
We also take the time
to explain the notions of phase diagram and
phase regions. We believe this paper will
be helpful to readers who want to understand HCT
and its performance, but who  would be overwhelmed
by the epsilontics of the analysis in \cite{AoS}.

{ The paper is organized as follows.    Section  \ref{sec-ideal}   introduces a functional framework and several ideal quantities.   These include  the proxy classification error where Fisher's separation (SEP) \cite{Anderson}  plays a key role,   the  ideal  threshold  as a proxy for the optimal threshold, and the ideal HCT  as a proxy of the HCT.         Section  \ref{sec-HCT}   introduces the main results on the  asymptotic behavior of the HC threshold 
under the asymptotic RW model,  and the focal point is the  phase diagram.         Section  \ref{sec-clipping}  outlines the basic idea behind the main results followed by the proofs. Section  \ref{sec-HCobjective}  discuss the connection between the ideal  threshold and the ideal HCT.     Section \ref{sec-FDR} discusses the ideal behavior of Bonferroni threshold feature selection
and FDR-controlling feature selection.    
Section \ref{sec-discussion} discusses 
the link between ideal HCT and ordinary HCT, the  
finite $p$ phase diagram, and  other 
appearances of HC in  recent literature.    }

\section{Sep functional and ideal threshold}  
\label{sec-ideal}
\setcounter{equation}{0} 
Suppose $L$ is a fixed, {\it nonrandom} linear classifier, 
with decision boundary $L <> 0$.
Will $L$ correctly classify the future realization $(Y,X)$
from simple model of Section \ref{sec-simple-model}.
Then $Y = \pm 1$ equiprobable and $X \sim N(Y\mu,I_p)$.
The misclassification probability can be written 
\begin{equation}  \label{Relationship1}
   P\{ YL(X) < 0  | \mu\} = \Phi(-  \frac{1}{2} Sep(L;\mu)),
\end{equation} 
where $\Phi$ denotes the standard normal distribution function and 
$ Sep(L;\mu)$ 
measures the standardized interclass
distance: 
\begin{eqnarray}  \label{DefineSep1} 
   Sep(L;\mu) &=& \frac{ E\{L(X) |   Y=1\} - E\{L(X) |  Y=-1\} }{SD(L(X))}\\
                &=& \frac{   2 \sum w(j) \mu(j)  }{(\sum w(j)^2)^{1/2}}  = \frac{ 2  \langle w, \mu\rangle  }{\|w\|_2},   \nonumber
\end{eqnarray} 
The {\it ideal linear classifier} $L_\mu$ 
with feature weights $w \propto \mu$, and decision threshold
$L_\mu <> 0$ implements the likelihood ratio test.  It also 
maximizes $Sep$, since for every other linear classifier $L$,
$Sep(L; \mu) \leq Sep(L_\mu; \mu) = 2 \|\mu\|_2$.

\subsection{Certainty-equivalent heuristic}
\newcommand{\proxySep}{\widetilde{Sep}}

Threshold selection rules  
give {\it random} linear classifiers:
the classifier weight vector $w$ is  a {\it random} variable,
because it depends on the $Z$-scores of the realized training sample. 
If $L_Z$ denotes a linear classifier constructed based 
on such a realized vector of $Z$-scores, then  the misclassification error
can be written as 
\begin{equation}  \label{Relationship2}
   Err(L_Z,\mu) =  \Phi(- \frac{1}{2} Sep(L_Z;\mu));
\end{equation} 
this is a random variable depending on $Z$ and on $\mu$.
 Heuristically, because there is a large number of
coordinates, some statistical regularity appears, and we anticipate that
random quantities can be replaced by expected values.
We proceed as if
\begin{equation} \label{SepHeuristic1}
    Sep(L_Z;\mu) \approx \frac{  2  E_{Z|\mu} \langle w, \mu\rangle  }{( E_{Z|\mu} w^2)^{1/2} }
\end{equation}
where the expectation is over the conditional distribution of $Z$ conditioned on $\mu$.
Our next step derives from the fact that $\mu$ itself is random,
having about $\eps \cdot p$ nonzero coordinates, in random locations.
Now as $w_i = \eta_t(Z_i)$ we write heuristically
\[
E_{Z|\mu} \langle w, \mu\rangle  \approx p \cdot \eps \cdot  \mu_0 \cdot E  \{ \eta_t(Z_1) | \mu_1=\mu_0 \},  
\]
while
\[
E_{Z|\mu} \langle w,w \rangle  \approx p \cdot  \left( \eps  \cdot E\{   \eta_t^2(Z_1) | \mu_1 = \mu_0 \}  + 
(1-\eps) \cdot   E\{   \eta_t^2(Z_1) | \mu_1 =0 \}  \right) .
\]

\begin{definition}
Let the threshold $t$ be fixed and chosen independently of the training set. 
In the $RW(\eps,\tau)$ model we use the 
following expressions for {\bf proxy separation} 
\[
      \widetilde{Sep}(t; \eps,\tau) = \frac{2A}{\sqrt{B}},
\]
where 
\[
    A(t,\eps,\tau)  =    \eps \cdot  \tau \cdot E  \{ \eta_t( \tau+ W ) \}  .
\]
\[
   B(t,\eps,\tau)  =    \eps \cdot E  \{ \eta_t^2( \tau+ W ) \}   + (1-\eps) E  \{ \eta_t^2(  W ) \} .
\]
and $W$ denotes a standard normal random variable.
By {\bf proxy classification error} we mean
\[
   \widetilde{Err}(t;\eps,\tau,p,n)= \Phi( -  \frac{1}{2} \sqrt{ \frac{p}{n}} \cdot  \widetilde{Sep}(t,\eps,\tau)).
\]
\end{definition}
Normalizations are chosen here so that, in large samples
\[
               Sep(L_Z,\mu) \approx \sqrt{\frac{p}{n}} \cdot \widetilde{Sep}(t,\eps,\tau)  . 
\]
While ordinarily, we expect averages to ``behave like expectations'' in large samples,
we use the word {\it proxy} to remind us that there is a difference (presumably small).
Software to compute these proxy expressions has been developed
by the authors, and some numerical results using them were
reported in \cite{PNAS}.

Of course, the rationale for our interest in these proxy expressions
is our heuristic understanding that they accurately describe the exact large-sample behavior 
of certain threshold selection schemes.  This issue is settled in the affirmative,
after considerable effort, in \cite{AoS}.

\subsection{Certainty-equivalent threshold functionals}

In general the best threshold to use in a given instance of the RW model
depends on both $\eps$ and $\tau$.  It also depends on the specific realization  
of $\mu$ and even of $Z$. However,  dependence on $\mu$ and $Z$ is simply ``noise''
that goes away in large samples, while the dependence on $\eps$ and $\tau$ remains.

\begin{definition}
The {\it ideal threshold functional}
$T_{ideal}(\eps,\tau)$ maximizes the proxy separation
\[
     T_{ideal}(\eps,\tau)  = \mbox{arg max}_t  \proxySep(t,\eps,\tau).
\]
\end{definition}
Heuristically, $T_{ideal}$ represents a near-optimal threshold
in all sufficiently large samples; it is what we ``ought'' to be attempting to
use.

\begin{definition} {\bf Folding.} \label{def-fold}
The following concepts and notations will be used in connection
with distributions of absolute values of random variables. 
\bitem
\item The Half Normal distribution function $\Psi(t) = P \{ |N(0,1)| \leq t \}$.
\item The noncentral Half-Normal distribution   $\Psi_\tau(t) = P \{ |N(\tau,1)| \leq t \}$.
\item Given a distribution function $F$, the {\it folded distribution} is $G(t) =  F(t) - F(-t)$.
The Half Normal is the folded version of the standard Normal, and
the noncentral Half Normal is the folded version of a Normal with unit standard deviation
and nonzero mean equal to the noncentrality parameter.
\item Let $F_{\eps,\tau}$ denote the 2-point mixture
\[
    F_{\eps,\tau} (t) = (1-\eps) \Phi(t) + \eps \Phi(t-\tau); 
\]
 $G_{\eps,\tau}$  denotes the corresponding folded distribution:
 \[
    G_{\eps,\tau} (t) = (1-\eps) \Psi_0(t) + \eps \Psi_\tau(t), 
\]
 \eitem
 \end{definition}

We now define an HCT functional
representing  the target that HC thresholding aims for.

\begin{definition}
Let $F$ be a distribution function which is not
the standard normal $\Phi$.
At such a distribution, we define the HCT functional by
\[
    T_{HC}(F)  = \mbox{argmax}_{t > t_0}  \frac{ \bar{G}(t) - \bar{\Psi}(t)}{\sqrt{G(t)\cdot \bar{G}(t)}} ;
\]
here $G$ is the folding of $F$,
and $t_0$ is a fixed parameter of the HCT method (eg.  $t_0 = \Phi^{-1}(0.1)$).
The HC threshold in the $RW(\eps,\tau)$ model 
may be written, in an abuse of notation, 
\[
    T_{HC}(\eps,\tau) 
\]
meaning $T_{HC}(F_{\eps,\tau})$.
\end{definition}

Let $F_{n,p}$ denote the usual empirical distribution
of the feature $Z$-scores $Z_i$. The HCT of Definition \ref{hct-def} can be written as
$\hat{t}_{n,p}^{HC} = T_{HC}(F_{n,p})$.  Let $F$ denote the expected value of 
$F_{n,p}$; then $T_{HC}(F)$ will be called the {\it ideal HC threshold}.
Heuristically, we expect the usual sampling fluctuations and that
\[
     T_{HC}(F_{n,p}) \approx T_{HC}(F)
\] 
with a discrepancy decaying as $p$ and $n$ increase. This issue
is carefully considered in the companion paper \cite{AoS}, which 
shows that  the empirical HC threshold
in the ARW model indeed closely matches  the ideal HC threshold.

For comparison purposes, we considered two other threshold schemes.
First, (ideal) False-Discovery Rate thresholding. For a threshold $t$,
and parameters $(p,\eps,\tau)$, the expected number of useful features selected is
\[
    E(TP) (t;\eps,\tau,p ) = p \cdot \eps \cdot \bar{G}_{\eps,\tau}(t) ;
\]
and the expected number of useless features selected is
\[
  E(FP) (t;\eps,\tau,p ) = p \cdot (1-\eps) \cdot \bar{\Psi}(t) .
\]
Let $TPR(t) = p^{-1} E(TP)(t)$ denote the
expected {\it rate} of useful features above threshold and
$FPR(t) = p^{-1} E(FP)(t)$ denote the expected
{\it rate} of usless features above threshold. 
In analogy with our earlier heuristic, we define the proxy False Discovery Rate (FDR)
\[
  \widetilde{FDR}(t;\eps,\tau,p) = \frac{FPR(t)}{TPR(t) + FPR(t)}
\]
(The term ``proxy'' reminds us that
\[
      \frac{E (FP)(t)}{E(TP)(t) + E(FP)(t)}  \neq E   \frac{FP(t)}{TP(t) + FP(t)},
\]
although for large $p$ the difference will often be small.)

We define the FDRT-$\alpha$ functional by
\[
  T_{FDR,\alpha}(\eps,\tau) = \min \{ t :  \widetilde{FDR} < \alpha,\quad  t > t_0 \} . 
\]
Heuristically, this is the threshold that FDRT is `trying' to learn
from noisy empirical data.
We will also need the {\it proxy Local FDR}.
\[
  \widetilde{Lfdr}(t;\eps,\tau,p) = \frac{FPR'(t)}{TPR'(t) + FPR'(t)}.
\]
Here $FPR'$ denotes the derivative of  $FPR$, which exists,
using smoothness properties of $\Psi_0$; similarly for $TPR'$ and $\Psi_\tau$.
Intuitively, $  \widetilde{Lfdr}(t)$ denotes the expected fraction of
useless features among those features having observed $Z$-scores near level $t$.

Second, we considered Bonferroni-based thresholding.
\[
 T_{Bon} = \bar{\Phi}^{-1}(p^{-1}).
\]
This threshold level is set at the level that would cause
on average one false alarm in a set of $p$ null cases.

In \cite[Figures 2-3]{PNAS}, we presented  
  numerical calculations
of all these functionals and their separation behavior in two cases.
\bitem
 \item   $p = 10,000$, and $\eps = .01$.
 \item   $p = 10^6$ and $\eps = .0001$.
\eitem 
Although our calculations are exact numerical finite-$p$ calculations,
we remark that they correspond to sparsity exponents
$\beta = 1/2$ and $\beta = 2/3$, respectively.   
The figures show the following.
\bitem
 \item There is a very close numerical approximation of the HCT to the
 ideal threshold, not just at large $\tau$ but also even at quite small
 $2 < \tau < 3$.
 \item FDR and Bonferroni thresholds behave very differently from
 the ideal and from HC.
 \item The separation behavior of the HCT is nearly ideal. For the
 constant FDR rules, the separation behavior is close to ideal at some
 $\tau$ but becomes noticeably sub-ideal at other $\tau$.
 \item The False discovery rate behavior of HCT and Ideal thresholding
 depends on $\tau$. At small $\tau$, both rules tolerate a high FDR
 while at large $\tau$, both rules obtain a small FDR.
 \item The Missed detection rate of HCT and Ideal thresholding 
 also depends on $\tau$. At small $\tau$, the missed detection rate
 is high, but noticeably less than 100\%. At large $\tau$, the missed
 detection rate falls, but remains noticeably above $0\%$. In contrast
 the MDR for FDR procedures is essentially 100\% for small $\tau$
 and falls below that of HCT/ideal for large $\tau$.
 \eitem
 
 These numerical examples illustrate the idealized behavior of different
 procedures.  We can think of the HCT functional as the threshold
 which is being estimated by the actual HCT rule.
On an actual dataset sampled from the underlying $F$, the HC threshold
will behave differently, primarly due to stochastic fluctuations $F_{n,p} \approx F$.
Nevertheless, the close approximation of the HCT threshold to the ideal
one is striking and, to us, compelling.

\section{Behavior of ideal threshold, asymptotic RW model}  \label{sec-HCT}
\setcounter{equation}{0} 
\newcommand{\proxyErr}{\widetilde{Err}}
We now study the ideal threshold
in the asymptotic RW model of Definition \ref{def-ARW}.  That is,
we fix parameters $r$,$\beta$ in that model and study  the 
choice of threshold $t$ maximizing class separation. 

We first make precise a structural fact  about the ideal
threshold, first observed informally in \cite{PNAS}.
\begin{definition} {\bf ROC Curve.}
The {\em feature detection receiver
operating characteristic curve} (ROC) is the curve parameterized
by $(FPR(t),TPR(t))$.  The tangent to this curve at $t$
is
\[
    tan(t) = \frac{TPR'(t)}{FPR'(t)}
\]
and the secant is 
\[
   sec(t) = \frac{TPR(t)}{FPR(t)}.
\]
Note that in the $RW(\eps,\tau)$ model,
 $TPR$, $FPR$, $tan$ and $sec$ all depend on
$t$, $\eps$,$\tau$ and $p$, although we may, as here,
indicate only dependence on $t$.
\end{definition}

 \begin{theorem} {\bf Tangent-Secant Rule.}
In the $ARW(r,\beta)$ model. we have
\begin{equation} \label{SecTanRule}
       \frac{tan(T_{Ideal})}{1 + tan(T_{Ideal})} 
        \sim \frac{1}{2}  \left( 1 +  \frac{sec(T_{Ideal})}{1 + sec(T_{Ideal})}  \right), \qquad p \goto \infty.
\end{equation}
Here  $\eps = p^{-\beta}$, $\tau = \sqrt{2r \log(p)}$ and  $n \sim c \log(p)^\gamma$, $p \goto \infty$ as
in Definition \ref{def-ARW}. 
\end{theorem}

\begin{definition} {\bf Success Region.}
The {\em{region of asymptotically successful ideal threshold feature selection}} 
in the $(\beta,r)$ plane is the interior of the subset where
the ideal threshold choice $T_{ideal}(\eps,\tau)$ obeys
\[
\proxyErr(T_{ideal}(\eps,\tau);\eps,\tau,p) \goto \infty,    \;\;\; p \goto \infty ;
\]
here  we are in the $ARW(r,\beta)$ model of Definition \ref{def-ARW}.
\end{definition}
The interesting range involves $(\beta,r) \in [0,1]^2$.
The following function is important for
our analysis, and has previously appeared
in important roles in other (seemingly unrelated) problems;
see Section \ref{sec-OtherHC}.
\begin{equation}    \label{eq-rho-def}
\rho^*(\beta) =  
\left\{
\begin{array}{ll}
0,  &\qquad 0 < \beta \leq  1/2,  \\
\beta - 1/2,  &\qquad 1/2 < \beta \leq  3/4,  \\
(1 - \sqrt{1 - \beta})^2,   & \qquad  3/4 < \beta < 1. 
\end{array} 
\right. 
\end{equation}
As it turns out, it marks the boundary between
success and failure for threshold feature selection.


\begin{theorem} {\bf Existence of Phases.}
The success region is precisely
$r  > \rho^*(\beta)$, $0 < \beta < 1$.  
In the interior of the complementary 
region $r < \rho^*(\beta)$, $1/2 < \beta < 1$, 
even the ideal threshold cannot
send the proxy separation to infinity with increasing $(n,p)$.
\end{theorem}

\begin{figure}[h]   
\centering
\includegraphics[height = 4 in, width = 5  in]{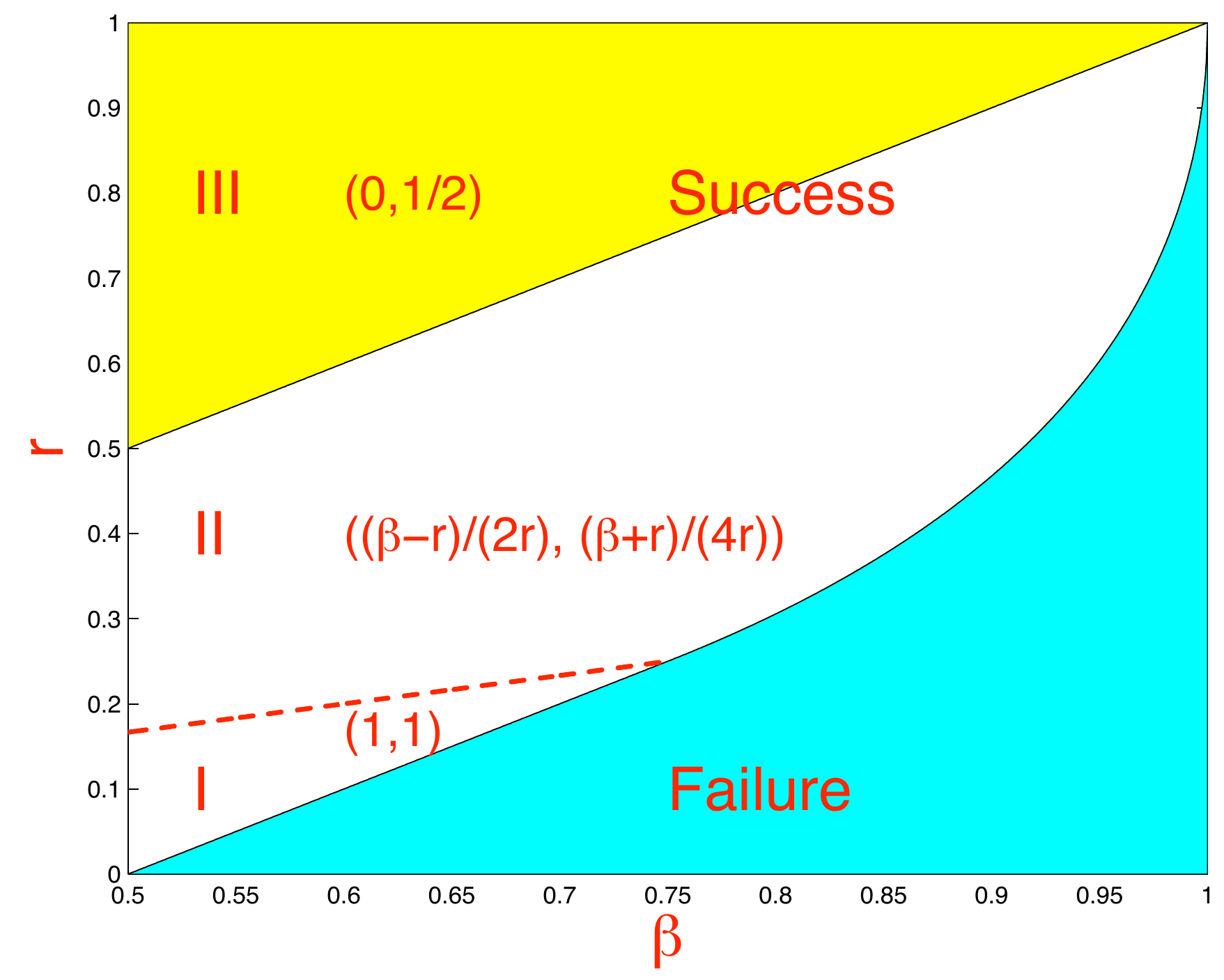}
\caption{Phase diagram.   The curve $r = \rho^*(\beta)$ splits the phase space into the failure region and the success region, and the latter further splits into three different regions I, II, III.    Numbers in the brackets show limits  of  
$\widetilde{FDR}$ and $\widetilde{Lfdr}$ at the ideal HCT as in Theorem 4. }
\label{fig:phase} 
\end{figure}

\begin{definition} \label{def-regions} {\bf Regions I,II, III.}
The  Success Region
can be split into three regions, referred to here and below as Regions I-III.
The interiors of the regions are as follows:
\begin{description}
\item[I.]   $\beta - 1/2 < r \leq  \beta/3$ and $1/2 < \beta < 3/4$; $r > \rho^*(\beta)$.  
\item[II.]  $\beta/3 <  r \leq  \beta$ and  $1/2 < \beta <1$; $r > \rho^*(\beta)$. 
\item[III.]   $\beta<  r  <   1$ and  $1/2 < \beta <1$; $r > \rho^*(\beta)$. 
\end{description} 
See Figure \ref{fig:phase}.  
\end{definition}
 
In the asymptotic RW model, the optimal threshold must behave asymptotically
like $\sqrt{2q \log(p)}$ for a certain $q = q(r,\beta)$.
Surprisingly we need not have $q = r$.
\begin{theorem}  {\bf Formula for Ideal Threshold.}
Under the Asymptotic RW model  $ARW(r,\beta)$, 
with $r > \rho^*(\beta)$,
the ideal threshold has the form $T_{ideal}(\eps,\tau)   \sim  \sqrt{2q^*\log(p)}$ where 
\begin{equation}  \label{qformula}
q^*   =   \left\{
\begin{array}{ll}
4r,  &\qquad        \mbox{Region I},    \\
\frac{(\beta + r)^2}{4r},  &\qquad \mbox{Region II,III}.       \\
\end{array}  
\right. 
\end{equation} 
\end{theorem}
Note in particular that in Regions I and II,
 $q^* > r$, and hence $T_{ideal}(\eps,\tau)  > \tau$.  Although the features truly have strength $\tau$,
the threshold is best set {\it higher than} $\tau$.

We now turn to FDR properties.
The Tangent-Secant rule implies immediately
\begin{equation} \label{eq-LFDR-FDR}
      \frac{Lfdr(T_{Ideal})}{(1 + FDR(T_{Ideal}))/2}    \goto 1, \qquad p \goto \infty.
\end{equation}
Hence any result about FDR is tied to one about  local FDR,
and vice versa.

\begin{theorem}
Under the Asymptotic RW model  $ARW(r,\beta)$,  
at the ideal threshold $T_{ideal}(\eps,\tau)$ proxy FDR obeys
\begin{equation} \label{Relationship5} 
\widetilde{FDR}(T_{ideal},\eps,\tau)    \goto
\left\{
\begin{array}{ll}
  1,  &\qquad    \mbox{Region I} \\
\frac{\beta - r}{2r},  &\qquad     \mbox{Region II}   \qquad  (\mbox{note $\beta/3 < r < \beta$})   \\
0,    &\qquad    \mbox{Region III}
\end{array} 
\right.  
\end{equation} 
as $p \goto \infty$, and the proxy local FDR 
obeys
\begin{equation} \label{Relationship6} 
\widetilde{Lfdr}(T_{ideal},\eps,\tau)  \goto 
\left\{
\begin{array}{ll}
  1,  &\qquad    \mbox{Region I} \\
\frac{ r + \beta }{4r},  &\qquad     \mbox{Region II}  \qquad    (\mbox{note $\beta/3 < r < \beta$})    \\
1/2,    &\qquad    \mbox{Region III}
\end{array} 
\right.  
\end{equation} 
as $p \goto \infty$.
\end{theorem}

Several aspects of the above solution are of interest.
\bitem
\item {\it Threshold Elevation.} The threshold $\sqrt{2q^* \log(p)}$ is significantly higher than
$\sqrt{2r\log(p)}$ in Regions I and II.  Instead of looking for features
at the amplitude they can be expected to have, we look for them at much higher 
amplitudes.
\item {\it Fractional Harvesting.} Outside of Region III,  we are selecting
only a small fraction of the truly useful features.
\item {\it False Discovery Rate.}  Outside Region III,   we actually have a very large false discovery rate, which is very close to $1$ in Region I.  Surprisingly {\it even though most of the selected features are useless},
we still correctly classify!   
\item {\it Training versus Test performance.}
The quantity $\sqrt{q}$ can be interpreted as a ratio: $\sqrt{q^*/r}  =  \min\{  \frac{\beta + r}{2r}, 2\} = $ strength of useful features in
training / strength of those features in test. From (\ref{qformula}) we learn that, in Region I, 
the selected useful features perform about half as well in training as we might
expect from their performance in test.

\eitem

\section{Behavior of  ideal clipping threshold} \label{sec-clipping} 
\setcounter{equation}{0} 

We now sketch some of the arguments involved
in the full proof of the theorems stated above.
In the RW model, it makes particular 
sense to use the clipping threshold function
$\eta_t^{clip}$, since all nonzeros 
are known to have the same amplitude.
The ideal clipping threshold is also very easy to analyze
heuristically.     But it turns out that all the statements in Theorems 1-3
are equally valid for all three types of threshold functions, so we prefer
to explain the derivations using clipping.

\newcommand{\sgn}{\mbox{sgn}}
\subsection{$\proxySep$ in terms of true and false discoveries}
In the RW model, we can express the components
of the proxy separation very simply when using the
clipping threshold:
\[
   A_{clip}(t,\eps,\tau) = \eps  \cdot \tau \cdot E \sgn(\tau + W) 1_{\{ |\tau + W | > t \}}.
\]
\[
   B_{clip}(t,\eps,\tau) = \eps\cdot E 1_{\{ |\tau + W | > t \}} + (1-\eps) \cdot E 1_{\{ | W | > t \}}
\]
where $W$ denotes an $N(0,1)$ random variable. Recall the definitions of useful selections TP
and useless selections FP; we also must count {\it Inverted Detections},  for the case where
the $\mu_i > 0$ but $\eta_t^{clip}(Z_i) < 0$.  Put
\[
   E(ID)(t;\eps,\tau,p) = \eps \cdot p \cdot \Phi(-t-\tau),
\]
with again $\Phi$ the standard normal distribution,
and define the {\it inverted detecton rate} 
by $IDR = p^{-1} E(ID)$.
Then
\[
    A_{clip}(t;\eps,\tau) = \tau\cdot \left(TPR(t) - 2 IDR)(t) \right) .
\]
\[
    B_{clip}(t;\eps,\tau) = TPR(t) + FPR(t) .
\]
We arrive at an identity for $\proxySep$ in the case of clipping:
\[
\proxySep(t ; \eps,\tau) = \frac{(2 \sqrt{p} \tau)  \cdot \left(TPR(t) - 2IDR(t) \right) }{\sqrt{ TPR(t) + FPR(t)}}.
\]

\newcommand{\proxyoSep}{\overline{Sep}}
We now explain Theorem 1, the Tangent-Secant rule.
Consider the alternate proxy
\[
   \proxyoSep(t) = \frac{ (2  \tau)  \cdot TPR(t)  }{\sqrt{ TPR(t) + FPR(t)}} =  \frac{2 A(t) }{B^{1/2}(t) },  \mbox{ say};
\]
i.e.  drop the term $IDR$.
It turns out that for the alternate proxy, the tangent secant rule and
resulting FDR-Lfdr balance equation are exact identites.

\begin{lemma}
Let $\eps > 0$ and $\tau > 0$.
The threshold  $t_{alt}$ maximizing 
$   \proxyoSep(t) $ as a function of $t$
satisfies the Tangent-Secant rule as an exact identity;
at this threshold we have
\begin{equation} \label{alt-fdr-lfdr}
    Lfdr(t_{alt}) = \frac{1}{2} \left( 1 + FDR(t_{alt}) \right).
\end{equation}
\end{lemma}

{\bf Proof.} 
Now, $A$ and $B$ are both smooth functions of $t$, so
at the $t$ optimizing $A B^{- 1/2}$ we have
\[
0  =  B^{-1/2} \left( A' - \frac{1}{2} \frac{A}{B} B' \right)    =  
    \frac{B'}{B^{1/2}} \cdot \left( \frac{A'}{B'} - \frac{1}{2} \frac{A}{B} \right). 
\]
By inspection $B'(t) < 0$ for every $t > 0$. Hence,
\begin{equation} \label{eq-AB-identity}
   \frac{A'}{\tau B'} = \frac{1}{2} \frac{A}{\tau B} .
\end{equation}
The Tangent-Secant Rule follows.
We now remark that
\[
    FDR(t) = 1 - \frac{ \ TPR(t)  }{{ TPR(t) + FPR(t)}}  = 1 - \frac{A}{\tau \cdot B} ,
\]
and
\[
    Lfdr(t) = 1 - \frac{ \ TPR'(t)  }{{ TPR'(t) + FPR'(t)}} = 1 - \frac{A'}{\tau \cdot B'}  .
\]
Display (\ref{alt-fdr-lfdr}) follows \qed

The full proof of Theorem 1, which we omit,
simply shows that the discrepancy caused by
$ \proxyoSep \neq  \proxySep$ has an asymptotically negligible
effect; the two objectives have very similar maximizers.

\subsection{Analysis in the asymptotic RW model}

We now invoke the $ARW(r,\beta)$ model: $\eps = p^{-\beta}$,
$\tau = \sqrt{2 r \log(p)}$, $(n,p) \goto \infty$, $n \sim c \log(p)^\gamma$, $p \goto \infty$.
Let $t_p(q) = \sqrt{2 q \log(p)}$. 
The classical Mills' ratio can be written
in terms of the Normal survival function as:
\[
     \bar{\Phi}(t_p(q)) \sim p^{-q} /t_p(q),   \qquad p \goto \infty.
\]
Correspondingly, the half Normal obeys:
\begin{equation} \label{mills}
     \bar{\Psi}_0(t_p(q)) \sim 2 \cdot p^{-q} /t_p(q),   \qquad p \goto \infty.
\end{equation}

\newcommand{\PL}{\mbox{\sc PL}}
We also need a notation for poly-log terms.
\begin{definition} \label{def-PL}
Any occurrence of the symbol $\PL(p)$ denotes a term 
which is $O(\log(p)^\zeta)$ and  $\Omega(\log(p)^{-\zeta})$  as $p \goto \infty$
for some $\zeta > 0$.  Different occurrences of this symbol may stand for different
such terms.
\end{definition}
In particular, we may well have $T_1(p) = \PL(p)$,  $T_2(p) = \PL(p)$,  as $p \goto \infty$,
and yet $\frac{T_1(p)}{T_2(p)} \not \goto 1$ as $p \goto \infty$.  
However, certainly $\frac{T_1(p)}{T_2(p)} = \PL(p)$, $p \goto \infty$.

The following Lemma exposes the main phenomena
driving  Theorems 1-4.  It follows by simple algebra,
and several uses of Mills' Ratio (\ref{mills})
in the convenient form $ \bar{\Psi}_0(t_p(q)) =  PL(p) \cdot p^{-\beta}$.
\begin{lemma} \label{MainLemma}
In the asymptotic RW model $ARW(r,\beta)$,
we have: 
\begin{enumerate}
\item {\it Quasi power-law for  useful feature discoveries}:
\begin{equation}
E(TP)(t_q(p),\eps,\tau) = \PL(p) \cdot p^{\delta(q,r,\beta)} , \qquad p \goto \infty.
\end{equation}
where the useful feature discovery exponent $\delta$ obeys
\[
\delta(q; \beta,r) \equiv   \left\{
\begin{array}{ll}
{1 - \beta},  &\qquad 0 < q \leq r,  \\
{1 - \beta - (\sqrt{q} - \sqrt{r})^2},   &\qquad  r <  q < 1. 
\end{array} 
\right. 
\]
\item {\it Quasi power-law for useless feature discoveries}:
\begin{equation}
E(FP)(t_q(p),\eps,\tau) = \PL(p) \cdot p^{1-q}, \qquad p \goto \infty,
\end{equation}
\item {\it Negligibility of inverted detections}:
\begin{equation}
E(ID)(t_q(p),\eps,\tau) = o( E(TP) (t_q(p),\eps,\tau) ) , \qquad p \goto \infty.
\end{equation}
\end{enumerate}
\end{lemma}

As an immediate corollary,  under  $ARW(r,\beta)$, we have:
\begin{equation}\label{eq-asysep}
\proxySep(t_q(p),\eps,\tau) \cdot (\frac{1}{2\sqrt{n}})    = \PL(p) \cdot   \frac{p^{\delta(q;\beta,r)}}{\sqrt{p^{\delta(q;\beta,r)} +  p^{1 - q}}} , \qquad p \goto \infty.
\end{equation}
On the right side of this display,
the poly-log term is relatively unimportant.
The driving effect is the power-law behavior of the fraction.
The following Lemma contains the core idea behind
the appearance of $\rho^*$ in Theorems 1 and 2,
and the distinction between Regions I and Region II,III.
\newcommand{\proxyHC}{\widetilde{HC}}
\begin{lemma} \label{lem-ratiolemma}
Let $(\beta,r) \in(\frac{1}{2},1)^2$.
Let $\gamma(q;r,\beta)$ denote the rate at which
\[
 \frac{p^{\delta(q;\beta,r)}}{\sqrt{p^{\delta(q;\beta,r)} +  p^{1 - q}}}  
\]
tends to $\infty$ as $p \goto \infty$, for fixed $q$, $r$, and $\beta$.
Then $\gamma > 0$ if and only if $r > \rho^*(\beta)$.  A choice of
$q$ maximizing this ratio is given by (\ref{qformula}).
\end{lemma}

{\bf Proof Sketch.}  By inspection of $\delta$, it is enough to consider $q \leq 1$.
The ratio grows to infinity with $p$ like $p^{\gamma}$, where 
\[
\gamma(q;r,\beta) =  \delta(q;\beta,r)  - \max((1-q),\delta(q;\beta,r))/2 ;
\]
provided there exists  $q \in [0,1]$  obeying $\gamma(q;r,\beta) > 0$. 

Let $q^*(r,\beta)$ denote the value of $q$ maximizing the rate of separation:
\[
   q^*(r,\beta) = \mbox{argmax}_{q \in [0,1]}  \gamma(q;r,\beta) ;
\]  
in the event of ties, we take the smallest value of $q$.
Let's define $\rho^*(\beta)$
without recourse to the earlier formula (\ref{eq-rho-def})  but instead by the functional 
role claimed for it by this lemma:
\begin{equation} \label{alt-rho-def}
  \rho^*(\beta) = \inf \{ r : \gamma(q^*(r,\beta);r,\beta) > 0, r > 0 \}.
\end{equation}
We will derive the earlier formula (\ref{eq-rho-def})  from this.
  Now $\gamma = \min(\gamma_1,\gamma_2)$ where
\[ 
 \gamma_1(q) =  \delta(q;\beta,r)  - (1-q)/2; \qquad \mbox{ and  } \gamma_2(q) =  \delta(q;\beta,r) /2 .
\]
We have
\begin{equation} \label{eq-max-rep}
  \gamma(q^*(r,\beta);r,\beta)  = \max_{q \in [0,1]} \min_{i=1,2} \gamma_i(q;r\beta) .
\end{equation}

In dealing with this maximin,
two special choices of $q$ will recur below.
\bitem
\item $q_1$: Viewed as a function of $q$, $\gamma_1$ 
is maximized on the interval $[r,1]$ (use calculus!) at $q_1(r,\beta) \equiv 4r$, 
and is monotone on either side of the maximum.

\item $q_2$: On the other hand, $\gamma_2$ is monotone decreasing
as a function of  $q$ on $[r,1]$.
Hence the maximizing value of $q$  in (\ref{eq-max-rep})  
{\it over the set of $q$-values where $ \gamma_2$ achieves the minimum}
will occur  at the minimal value of $q$ achieving the minimum, i.e.
at the solution to:
\begin{equation} \label{eq-qdef}
   (1-q) = \delta(q;\beta,r) .
\end{equation}
(\ref{eq-qdef}) is satisfied uniquely 
on $[r,1]$ by $q_2(r,\beta) = (\beta+r)^2/4r$.
\eitem

The behavior of $\min(\gamma_1 ,\gamma_2)$ varies by cases;
see Table \ref{table-maximin}. To see how the table was derived,
note that
\[
    \gamma_1  < \gamma_2 \mbox{ iff } \delta < 1-q.
\]
Consider the first and second rows. For $q > r$, 
$\delta = 1 - q - \beta -r + 2\sqrt{rq}$. Hence
$\delta < 1-q$ on $[r,1]$  iff  $- \beta -r + 2\sqrt{rq} < 0$ iff  $ q < q_2$.
Consider the third and fourth rows.  For $q < r$,
 $\delta  = 1-\beta$. Hence $\delta < 1-q$ on $[0,r]$  iff $\beta > q$.
\begin{table}[h]
\begin{center}
\begin{tabular}{|l|l|}
\hline
Range & Minimizing $\gamma_i$\\
\hline
 $r \leq q \leq q_2$ & $\gamma_1$ \\
 $\max(r,q_2) < q $ & $\gamma_2$ \\
$\beta < q < r $ & $\gamma_2$ \\
 $q < \min(r,\beta)  $ & $\gamma_1$ \\
\hline
\end{tabular}
\caption{Minimizing $\gamma_i$ in (\ref{eq-max-rep})}
\end{center}
\label{table-maximin}
\end{table}

Derivatives 
$\dot{\gamma}_i = \frac{\partial}{\partial q} \gamma_i$, $i =1,2$,
are laid out in Table \ref{table-gammadot}.

\begin{table}[h]
\begin{center}
\begin{tabular}{|l|l|l|}
\hline
                                 &   $q < r$& $q > r$ \\
\hline
$\dot{\gamma}_1$ &$\frac{1}{2}$ & $-\frac{1}{2} + \sqrt{r/q}$ \\
$\dot{\gamma}_2$ &  0                & $-\frac{1}{2} + \frac{1}{2} \sqrt{r/q}$ \\
\hline
\end{tabular}
\caption{Derivatives of $\gamma_i$}
\end{center}
\label{table-gammadot}
\end{table}

Table \ref{table-maxq}
presents results of formally combining the two previous tables.
There are four different
cases, depending on the ordering
of $q_1$,$q_2$,$\beta$ and $r$.
In only one case does the above information leave
 $q^*$ undefined.  (We note that this is a purely formal calculation;
Lemma \ref{lem-algebra}, Display (\ref{eq-r-q2})
 below shows that {\it rows 2 and 4 never occur}.)
 
 To see how  Table \ref{table-maxq} is derived,
 consider the first row. Using the derivative table above,
 we see that $\min(\gamma_1,\gamma_2)$ is increasing
 on $[0,\beta]$, constant on $[\beta,r]$  and decreasing on $[r,1]$. Hence
 the maximin value is achieved at any $q \in [\beta,r]$.
For row  2, 
$\min(\gamma_1,\gamma_2)$ is increasing on $[0,\beta]$, constant on 
 on $[\beta,r]$ 
 increasing on $[r,\min(q_1,q_2)]$
 and decreasing on $[\min(q_1,q_2),1]$.
 For row  3, 
$\min(\gamma_1,\gamma_2)$ is constant on $[0,r]$, 
increasing on $[r,\min(q_1,q_2)]$
 and monotone decreasing on $[\min(q_1,q_2),1]$.
For row 4, $\min(\gamma_1,\gamma_2)$
 is increasing on $[0,r]$ and decreasing on $[r,1]$.

\begin{table}[h]
\begin{center}
\begin{tabular}{|l|l|l|l|}
\hline
Case & Minimizing $\gamma_i$   & Maximin $q$     & Maximin value \\
\hline
$\beta < r$, $q_2 < r$ & $ \left\{ \begin{array}{ll} 
                                                           \gamma_1 &  (0,\beta) \\
                                                           \gamma_2&  (\beta,1) \\
                                                           \end{array}  \right. $   
                                        & $q^* \in[\beta,r]$   & $\gamma_1(\beta), \gamma_2(q_2)$ \\
\hline                                                           
$\beta < r$, $q_2 > r$ & $ \left\{ \begin{array}{ll} 
                                                           \gamma_1 &  (0,\beta) \\
                                                           \gamma_2&  (\beta,r) \\
                                                            \gamma_1&  (r,q_2) \\
                                                            \gamma_2&  (q_2,1) \\
                                                           \end{array}  \right. $ 
                                        & $q^* =  \min(q_1,q_2)$ & $\gamma_1(q^*)$ \\
\hline 
$\beta > r$, $q_2 > r$ & $ \left\{ \begin{array}{ll} 
                                                           \gamma_2 &  (0,r) \\
                                                           \gamma_1&  (r,q_2) \\
                                                            \gamma_2&  (q_2,1) \\
                                                           \end{array}  \right. $
                                        & $q^* =  \min(q_1,q_2)$ & $\gamma_1(q^*)$ \\
\hline
$\beta > r$, $q_2 < r$ & $ \left\{ \begin{array}{ll} 
                                                           \gamma_1&  (0,r) \\
                                                           \gamma_2 &  (r,1) \\
                                                           \end{array}  \right. $
                                        & $q^* =r$ & $\gamma_1(r)$ \\
 \hline
\end{tabular}
\caption{Maximin Behavior in (\ref{eq-max-rep})}
\end{center}
\label{table-maxq}
\end{table}


We are trying to find circumstances
where $\gamma \leq 0$.
In the above table, we remarked that
the hypotheses of rows 2 and 4 can never occur.
We can see that in row 1, $\gamma_1(\beta; r,\beta) > 0$ for $\beta \in [0,1]$, $r > \beta$.
This leaves only row 3 where we might have $\gamma \leq 0$; in that case
either $q^* = q_1$ or $q^*=q_2$.
Writing out explicitly
\[
  \gamma_1(q;r,\beta) = 1 - \beta  +  ( \sqrt{q} - \sqrt{r})_+^2 - (1-q)/2 .
\]
\bitem
\item Case $q^* = q_1$.
\[
\gamma_1(q_1;r,\beta) = 1/2 - \beta  + r .
\]
Hence $\gamma_1(q_1;r,\beta)  = 0$
along $r = \beta - 1/2$, and $\gamma_1(q_1;r,\beta) < 0$
for $r < \beta -1/2$. 

Consider $1/2 < \beta < 3/4$. In this range,
Lemma \ref{lem-algebra} shows  $r < q_1( \beta -1/2,\beta) < q_2(\beta-1/2,\beta) < 1$.
Hence $q^*(r,\beta) = q_1(r,\beta)$  for $r \leq \beta - 1/2$, $\beta \in (1/2,3/4)$.
We conclude that
\begin{equation} \label{rho-formula-1}
  \rho^*(\beta) = \beta - 1/2, \qquad 1/2 \leq \beta \leq 3/4.
\end{equation}

\item Case $q^* = q_2$.
\[
\gamma_1(q_2;r,\beta) = 1/2  - (\beta+r)^2/8r.
\]
We have $\gamma_1(q_2;r,\beta)  = 0$
along $r = (1 - \sqrt{1-\beta})^2 $, and $\gamma_1(q_2;r,\beta) < 0$
for $0 < r < (1 - \sqrt{1-\beta})^2 $. 

Consider $3/4 < \beta < 1$. In this range,
Lemma \ref{lem-algebra} shows   shows $r < q_2((1 - \sqrt{1-\beta})^2,\beta) < q_1((1 - \sqrt{1-\beta})^2,\beta) < 1$.
We conclude that
\begin{equation} \label{rho-formula-2}
  \rho^*(\beta) = (1 - \sqrt{1-\beta})^2, \qquad 3/4 \leq \beta \leq 1/2.
\end{equation}

\eitem
Together (\ref{alt-rho-def}), (\ref{rho-formula-1}), and (\ref{rho-formula-2}) establish 
that the formula (\ref{eq-rho-def}) has the properties implied by the Lemma.

To complete the Lemma, we need to validate formula (\ref{qformula}).
This follows from rows 1 and 3 of
Table \ref{table-maxq} and Lemma \ref{lem-algebra}.
\qed

\begin{lemma} \label{lem-algebra}
Let $q_1(r,\beta) = 4r$ and $q_2(r,\beta) = (\beta+r)^2/4r$ just as in the previous lemma.
For $0 < \beta < 1$ we have
\begin{eqnarray}
q_1 < q_2 & \mbox{iff}  & r < \beta /3                      \nonumber \\
q_2 < 1     & \mbox{iff}  & r > (1 - \sqrt{1-\beta})^2 \nonumber \\
 r      < q_2& \mbox{iff}  & r < \beta                           \label{eq-r-q2}\\
\end{eqnarray}
\end{lemma}

{\bf Proof.}  School algebra. \qed

Lemmas \ref{lem-ratiolemma} and  \ref{lem-algebra} 
show that (\ref{qformula}) gives us {\it one choice} of  threshold
maximizing the rate at which $\proxySep$ tends to infinity;
is it the {\it only choice}? Except in the case,  $r > \beta$ and $q_2 < r$, 
this is indeed the only choice.
As  Table \ref{table-maxq} shows, in case  $r > \beta$ and $q_2 < r$, 
{\it any} $q \in [\beta,r]$ optimizes the {\it } of separation.
It turns out that in that case,  our formula $q^*$ not only
maximizes the rate of separation, it correctly describes the leading-order asymptotic behavior
of $T_{Ideal}$.  The key point is the Tangent-Secant Formula,
which picks out from among all $q \in [\beta,r]$ uniquely $q_2$.
This is shown by the next two lemmas, which thereby complete the proof of
Theorem 3.

\begin{lemma}
Set $\gamma_0(q;r,\beta) = -\beta - r + 2 \sqrt{rq}  $, for $q \in (0,1)$.
In the $ARW(r,\beta)$ model,
consider the threshold $t_q(p)$.
Suppose that $q > r$. Then 
\begin{equation} \label{eq-TP2FP}
    \frac{TPR(t_q(p))}{FPR(t_q(p))} \sim  p^{\gamma_0(q,r,\beta)} \cdot \frac{\sqrt{q}}{2\sqrt{q}-2\sqrt{r}} \qquad   p \goto \infty.
\end{equation}
Suppose that  $q <r$. Then
\begin{equation} \label{eq-TP2FPa}
    \frac{TPR(t_q(p))}{FPR(t_q(p))} \sim  p^{q-\beta} \cdot \sqrt{\frac{\pi}{2}} \cdot t_q(p)  \qquad   p \goto \infty.
\end{equation} 
Suppose that  $q \neq r$. Then
\begin{equation} \label{eq-LTP2LFP}
    \frac{TPR'(t_q(p))}{FPR'(t_q(p))}  =  \frac{1}{2} \cdot p^{\gamma_0(q,r,\beta)} .
\end{equation}
\end{lemma}

\noindent
{\bf Proof.}
Simple manipulations with Mills' ratio, this time {\it not} simply
grouping polylog terms together with the symbol $PL$, give that 
for the threshold under  the ARW model,  if $q > 0$,
\[
      FPR(t_q(p)) \sim\sqrt{\frac{2}{\pi}} \cdot    \frac{p^{-q}}{ \sqrt{2q \log(p)}} \qquad   p \goto \infty.
\]
If $q > r$
\[
    TPR(t_q(p)) \sim \sqrt{\frac{1}{2\pi}} \cdot \frac{p^{-\beta - (\sqrt{q}-\sqrt{r})^2}}{ (\sqrt{q}  - \sqrt{r}) \sqrt{ 2\log(p)}}.
\]
Display (\ref{eq-TP2FP}) follows.
If $0 < q < r$
\[
    TPR(t_q(p))) \sim \eps, \qquad   p \goto \infty.
\]
Display (\ref{eq-TP2FPa}) follows. If $q \neq r$ we have the exact identities:
\[
    TPR'(t_q(p)) = \sqrt{\frac{1}{2\pi}} \cdot p^{-\beta - (\sqrt{q} - \sqrt{r})^2 } , \qquad
      FPR'(t_q(p)) = \sqrt{\frac{2}{\pi}} \cdot    p^{-q}.
\]
Display (\ref{eq-LTP2LFP}) follows.
\qed

\begin{lemma}
$q_2(r,\beta)$ is the unique solution of $\gamma_0(q;r,\beta) = 0$.
Suppose $r > \beta$ and $\beta < q_2 < r$.
\begin{equation} \label{eq-ideal-lim}
   T_{Ideal} / t_{q_2}(p) \goto 1, \qquad p \goto \infty.
\end{equation}
\end{lemma}

\noindent
{\bf Proof.}
By inspection for $q < q_2$, $\gamma_0 < 0$
while for $q > q_2$, $\gamma_0 > 0$. So from (\ref{eq-LTP2LFP}),
$Lfdr(t_q(p))$ tends to $0$ or $1$ depending on $q < q_2$ or $q > q_2$.
Fix $\eta > 0$. If $T_{Ideal} < t_{q_2 - \eta}$ infinitely often
as $p \goto \infty$, then $0$ would be a cluster
point of  $Lfdr(T_{Ideal})$. Similarly, if $T_{Ideal} >t_{q_2 + \eta}$
infinitely often as $p \goto \infty$, then $1$ would be a cluster
point of  $Lfdr(T_{Ideal})$. From  $q > \beta$
and (\ref{eq-TP2FPa})  we know that  $FDR(T_{Ideal}) \goto 0$.
>From the Tangent-Secant rule we know that $Lfdr(T_{Ideal}) \goto 1/2$. 
(\ref{eq-ideal-lim}) follows.
\qed

\subsection{FDR/Lfdr  properties of ideal threshold}

We now turn to Theorem 4.
Important observation:
\begin{quotation}
{\it The asymptotic $T_{Ideal} =  t_{q^*}(p)\cdot (1+o(1))$ is simply too crude to determine 
the FDR and Lfdr properties of $T_{Ideal}$; it is necessary to consider
the second-order effects implicit in the $(1+o(1))$ term. For this, the
Tangent-Secant formula is essential.}
\end{quotation}

Indeed (\ref{eq-LTP2LFP}) shows that the only possibilities for
limiting local FDR of a threshold of the {\it exact} form $t_q(p)$
are $0, 1/2, 1$. The actual local FDR of $T_{Ideal}$ spans
a continuum from $[0,1]$, due to the fact that small perturbations
of $t_q(p)(1+o(1))$ implicit in the $o(1)$
can cause a change in the local FDR.
\newcommand{\tq}{\tilde{q}}
To understand this, for  $q \neq r$ and $s \in (0,\infty)$, put
\[
\tq(q,r,s,p) = \left( r^{1/2}  \pm \sqrt{ (q^{1/2} - r^{1/2})^2 + \log(s)/\log(p)} \right)^2 ,
\]
where the sign of $\pm$ is $+$ if $q > r$.  Clearly, $\tq$ is well-defined
for all sufficiently large $p$.
As an example, $\tq(q,r,1,p) = q$.  
The peculiar definition ensures that
\[
   \phi(t_{\tq}(p) - t_r(p)) = \phi(t_q(p)-t_r(p)) \cdot s.
\]

By simple algebra one can show
\begin{lemma}
Let $q$,$r$, and $s$ be fixed. With $\tq = \tq(q,r,s,p)$,
\[
   \frac{ t_{\tq}(p)}{t_q(p)} \goto 1, \qquad p \goto \infty.
\]
and
\[
   \frac{\phi( t_{\tq}(p))}{\phi(t_q(p))} \goto 1, \qquad p \goto \infty.
\]
\end{lemma}

Let's put for short $F_s = FPR(t_{\tq(q,r,s,p)}(p))$ and $T_s = TPR(t_{\tq(q,r,s,p)}(p))$.
Then, the last few displays show
\[
    T'_s = s \cdot T'_1 , \qquad F'_s \sim F'_1, \qquad  p \goto \infty.
\]
\[
    T_s \sim s \cdot T_1 , \qquad F_s \sim F_1, \qquad  p \goto \infty.
\]
Hence
\[
   FDR(t_{\tq(q,r,s,p)}(p)) = \frac{F_s}{T_s + F_s} \sim \frac{F_1}{s \cdot T_1 + F_1} , \qquad  p \goto \infty ;
\]
\[
   Lfdr(t_{\tq(q,r,s,p)}(p)) = \frac{F'_s}{T'_s + F'_s} \sim \frac{F'_1}{s \cdot T'_1 + F'_1}, \qquad  p \goto \infty.
\]

Choosing $s$ appropriately, we can therefore obtain a perturbed $q$
which perturbs the Lfdr and FDR. 
In fact there is a unique choice of $s$ needed to
ensure the Tangent-Secant formula.
\begin{lemma}
For given values $F_1$, $F'_1$, $T_1 \neq 0 $, $T'_1 \neq 0$, put
\[
    s^* = \frac{F'_1}{T'_1} - 2  \frac{F_1}{T_1}.
\]
This choice of $s$ obeys the Tangent-Secant rule:
\[
  \frac{F'_1}{s^* \cdot T'_1 + F'_1}  = \frac{1}{2} \left( 1 +   \frac{F_1}{s^* \cdot T_1 + F_1} \right).
\]
\end{lemma}

To use this recall (\ref{eq-TP2FP})-(\ref{eq-TP2FPa})-(\ref{eq-LTP2LFP}).
These formulas give expressions for $T_1/F_1$ and $T'_1/F'_1$.
Plugging in $q=q^*$ we get
\begin{corollary}
We have
\[
T_{Ideal} \sim t_{\tq(q^*,r,s^*,p)}(p) .
\]
where $s^*$ is obtained by setting $T_1 = (\sqrt{q^*}-\sqrt{r})^{-1}$, $F_1 = 2/\sqrt{q}$,
$T'_1=1$ and $F'_1= 2$. Moreover, if $\beta/3 < r < \beta$,
\[
    FDR(T_{Ideal}) \sim  \frac{\beta-r}{2r}, \qquad Lfdr(T_{Ideal}) \sim \frac{r+\beta}{4r} .
\]
\end{corollary}



\section{Connection of HC objective with  $\proxySep$}  \label{sec-HCobjective} 
\setcounter{equation}{0} 
\newcommand{\bPsi}{{\bar{\Psi}}}
Let $F = F_{\eps,\tau}$ be the two-point mixture of Definition \ref{def-fold} and 
$G = G_{\eps,\tau}$ the corresponding folded distribution. Then
in the asymptotic RW model we have, for $t = t_q(p)$:
\begin{eqnarray}
HC(t; F_{\eps,\tau})  &=& \frac{G(t) - \Psi(t)}{\sqrt{G(t)\bar{G}(t)}} \nonumber  \\
           &=& \frac{\eps (\Psi_\tau - \Psi_0)(t)}{[((1-\eps)\Psi_0 + \eps \Psi_\tau)((1-\eps)\bPsi_0 + \eps \bPsi_\tau)]^{1/2}} \nonumber \\
           &\sim& \frac{\eps (\Psi_\tau - \Psi_0)(t)}{[((1-\eps)\bPsi_0 + \eps \bPsi_\tau)]^{1/2}} \label{sim1} \\
           & \propto & \frac{TPR(t) - \frac{\eps}{1-\eps} FPR(t) }{[FPR(t)+ TPR(t) ]^{1/2}} \nonumber  \\
           & \sim &  \frac{TPR(t) }{[ TPR(t) + FPR(t)]^{1/2}} \label{sim2} \\
           &\sim & \proxySep(t;\eps,\tau), \qquad p \goto \infty .  \label{sim3}
\end{eqnarray}
Step (\ref{sim1}) follows from $G_{\eps,\tau}(t_q(p)) \goto 1$ as $p \goto \infty$.
Step (\ref{sim2}) follows from $\eps FPR(t_q(p))  = o(TPR(t_q(p)))$ as $p \goto \infty$.
Step (\ref{sim3}) follows from $TPR(t_q(p))  = o(TPR(t_q(p)))$ as $p \goto \infty$.

This derivation can be made rigorously correct when $q = q^*(r,\beta)$,
where $q^*$ is as announced in Theorem 2.  With extra work, not shown here,
one obtains:

\begin{theorem}
The statements made for the ideal threshold in Theorems 1-3
are equally valid for the ideal HC threshold.
\end{theorem}

\section{Suboptimality of phase diagram for other methods}  \label{sec-FDR}
\setcounter{equation}{0} 
Setting thresholding by control of the False Discovery rate
is a popular approach.  Another approach, more conservative
classical and probably even more popular,
is the Bonferroni approach, which controls the 
expected total number of false features selected.
Theorem 2 and 3 implicitly show the suboptimality of these approaches.
The implications include  an attractive relationship to the regions
of Definition \ref{def-regions}.

\begin{theorem}
Under the Asymptotic RW model  $ARW(r,\beta)$, 
the regions of the $(r,\beta)$ phase space
for successful classification are:
\[
r > (1 - \sqrt{1 - \beta})^2, \qquad  0 < \beta < 1.  
\]
\end{theorem}
{\bf  Remark.}  The successful region of Ideal FDRT and Bonferroni   is smaller than that of HCT.       However, even for regions when FDRT or Bonferroni would be successful, HCT stil yields more accurate classifications in terms of the  convergence rate. This is especially important for finite-p performance.    In short,   Bonferroni and FDRT fail  to adapt the difficulty level of the classification problem measured by $(\eps, \tau)$, one picks a fixed threshold, another picks a
fixed false discovery  rate.      

{\bf Proof Sketch.}
We continue to use the notations  $q^*(r,\beta)$, $q_i(r,\beta), i=1,2$, $\gamma(r,\beta)$
$\gamma_i(r,\beta)$, $i=1,2$ from the proof of Lemma \ref{MainLemma}.
The Bonferroni threshold $-\Phi^{-1}(1/p) = t_1(p)(1 + o(1))$  as $p \goto \infty$.  
In this proof sketch, we analyze $t_1(p)$ {\it as if} it were the Bonferroni threshold.
 Applying (\ref{eq-asysep}),
and the definition of $\gamma$ we have
\[
    \proxySep(t_1(p),\eps,\tau)  \cdot   \sqrt{\frac{p}{n}} 
    = PL(p)  \cdot p^{\gamma(1;r,\beta) }, \qquad p \goto \infty.
\]
Now 
\[
  \gamma(1;r,\beta) = \min_{i=1,2} \gamma_i (1 ; r, \beta ); 
\]
while $\gamma_1(1;r,\beta) = \delta(1;r,\beta)/2$ and $ \gamma_2(1;r,\beta) = \delta(1;r,\beta)$.
Hence $\gamma(1;r,\beta) > 0 $ iff $ \delta(1;r,\beta) > 0 $. But
\[
  \delta(1 ; r, \beta ) = -\beta -r + 2 \sqrt{r},
\]
which is positive iff $r > 1 - \sqrt{1-\beta}$.

Precise FDR control at level $\alpha$ requires that
$TP/FP = (1-\alpha)/\alpha$.  
Suppose that $r < \beta$.
Equation (\ref{eq-TP2FP}) relates the TP/FP
 ratio to $\gamma_0$, $q$ and $r$.   
 Clearly, the only way to keep
 $TP/FP$ bounded away from 0 and infinity is to choose $q$ so that
 $\gamma_0(q;r,\beta) = 0$.  We note that $q_2(r,\beta)$ exactly solves this problem:
 \[
   \gamma_0(q_2(r,\beta) ;r,\beta) = 0.
 \]
It follows that  in the $ARW(r,\beta)$ model,
the FDRT functional obeys
$T_{FDR,\alpha}(\eps,\tau) = t_{q_2}(p)(1 + o(1))$, $p \goto \infty$.
In this proof sketch we analyze $ t_{q_2}(p)$ {\it as if} it were 
exactly the FDR threshold.
\[
  \gamma(q_2(r,\beta);r,\beta) = \min_{i=1,2} \gamma_i (q_2(r,\beta) ; r, \beta ); 
\]
while $\gamma_1(q_2;r,\beta) =  1 - \frac{(\beta+r)^2}{4r}$ and $ \gamma_2(q_2;r,\beta) = \frac{1}{2} - \frac{(\beta+r)^2}{8r}$.
Hence $\gamma(q_2;r,\beta) > 0 $ is positive iff $r > 1 - \sqrt{1-\beta}$. 

The last paragraph assumed $r < \beta$. 
On inuitive grounds, the region $r > \beta$ offers even better
perfomance, so it must lie entirely above the phase transition. 
We omit details. \qed



\begin{figure} 
\begin{centering}
\includegraphics[height=3.5 in,width=4 in]{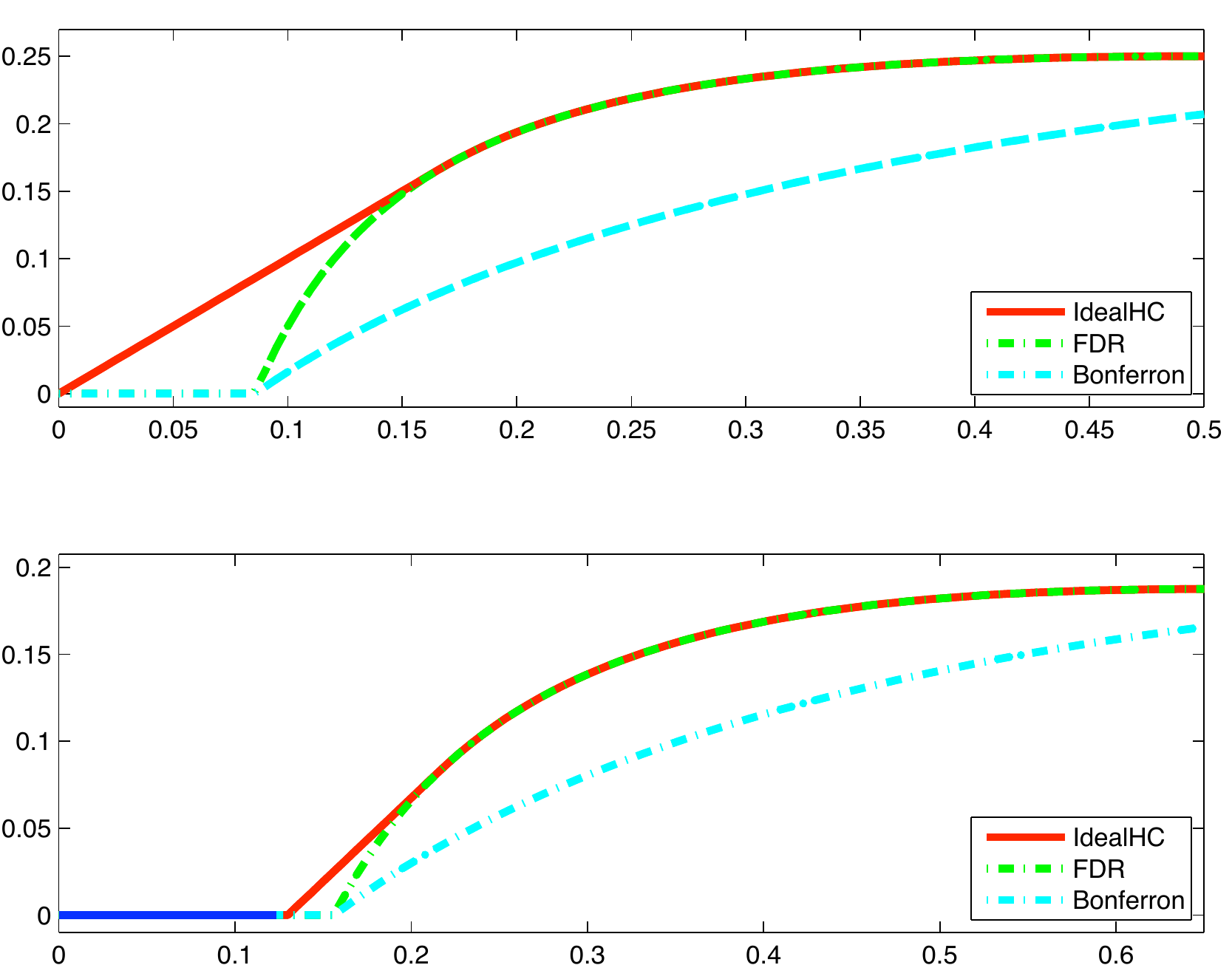}
 \caption{Comparison of SEP exponents for HCT,  FDRT, and Bonferroni.   
Top: $\beta = 1/2$. Bottom: $\beta =   5/8$.     
   $x$-axes displays $r$-values from $0$ to $\beta$, $y$-axes display corresponding exponents. Solid blue horizontal bar indicates that $r$ falls in  the failure region where all  exponents are $0$.  }
\label{fig:exponents}
\end{centering}
\end{figure}

Table \ref{table:exponents} compare the exponents in SEP for different methods. Also see Figure \ref{fig:exponents} for a comparison of  the exponents  for $\beta  = 1/2$ and $\beta = 5/8$.   

\begin{table}[htb]  
\begin{center}
\begin{tabular}{|r|l|l|l|}
\hline
Method   &   SEP Exponent & Boundary, $3/4 < \beta$ & Boundary, $1/2 < \beta < 3/4$ \\
\hline
Ideal, HCT        &   $\gamma$     &    $1 - \sqrt{1-\beta}$ & $\beta - 1/2$   \\
FDRT                &  $\frac{1}{2} - \frac{(\beta+r)^2}{8r} $  &  $1 - \sqrt{1-\beta}$ & $1 - \sqrt{1-\beta}$  \\
Bonferroni        & $ -\beta -r + 2 \sqrt{r}$ & $1 - \sqrt{1-\beta}$  &$1 - \sqrt{1-\beta}$  \\
\hline
\end{tabular}
\end{center}
\caption{Comparison of the exponents in SEP.} 
\label{table:exponents} 
\end{table}

\section{Discussion and conclusions}  \label{sec-discussion} 
\setcounter{equation}{0} 
\subsection{Beyond ideal performance}

We  consider here only asymptotic, ideal behavior.
Conceptually, the ideal threshold envisions a situation with
an oracle, who, knowing $\eps$ and $\tau$
and $n$ and $p$, chooses the very best threshold possible under those given parameters.
In this paper we have analyzed the behavior of this threshold within a
certain asymptotic framework. However, clearly, no empirical procedure
can duplicate the performance of the ideal threshold.

We have seen that ideal HC thresholding comes close. This ideal HC threshold does
not involve optimal exploitation of knowledge of $\eps$ and $\tau$,
but merely the availability of the underlying distribution of feature $Z$-scores
$F_{\eps,\tau}$. In this paper,we analyzed the behavior of $t^{HC} = T_{HC}(F_{\eps,\tau})$.

This is an ideal procedure because we never would know $F_{\eps,\tau}$;
instead we would have the empirical CDF $F_{n,p}$ defined by
\[
   F_{n,p}(z) = \frac{1}{p} \sum_{j=1}^p 1_{\{ Z(j) \leq z \} } .
\]
The (non-ideal) HCT that we defined in Section \ref{ssec-HCTDEF} is then simply
\[
    \hat{t}^{HC}_{n,p} = T_{HC}(F_{n,p}) .
\]
Because $F_{\eps,\tau}(z) = E(F_{n,p})(z)$, 
we are conditioned to expect that $t^{HC} \approx \hat{t}^{HC}_{n,p}$;
indeed there are generations of experience for {\it other} functionals
$T$ showing that we typically have $T(F_{n,p}) \approx T(E(F_{n,p}))$ for large $n$,$p$
for those functionals.
Proving this for $T = T_{HC}$ is more challenging than
one might anticipate; the problem is that $T_{HC}$ is
not continuous at $F = \Phi$, and yet $F_{\eps(p),\tau(p)} \goto \Phi$ as $p \goto \infty$.
After considerable effort, we justify the approximation of HCT by  ideal HCT 
in \cite{AoS}.  Hence the analysis presented here only partially proves
that HCT gives near-optimal threshold feature selection; it explains the
connection at the level of ideal quantities but not at the level of
fluctuations in random samples.

\subsection{Other asymptotic settings}

In the analysis here we consider only the case that $n \sim c \log(p)^\gamma$.
The phase diagram will be slightly different in case $n$ is bounded and
does not go to infinity with $p$, and will be again slightly different in case
$n \sim c p^\gamma$.  The full details are presented in \cite{AoS}.

\subsection{Phase diagram for finite sample sizes} 
{ While the main focus of our paper has been asymptotic analysis,   
we mention that the phase diagram also reflects  the finite-sample behavior.   
In Figure \ref{fig:finitep}, we consider    
$p = 3 \times 10^3 \times N$, $N = (1, 10, 100)$. 
For such $p$, we take $n = \log(p)/2$ and display the boundary of the set  of $(\beta, r)$ 
where ideal HCT yields a classification error between $10\%$ and $40\%$.   
The figure illustrates that as $p$ grows, both the upper bound and 
the lower bound migrate towards the common limit curve $r = \rho^*(\beta)$.

\begin{figure}[h]   
\centering
\includegraphics[height =4 in, width = 5  in]{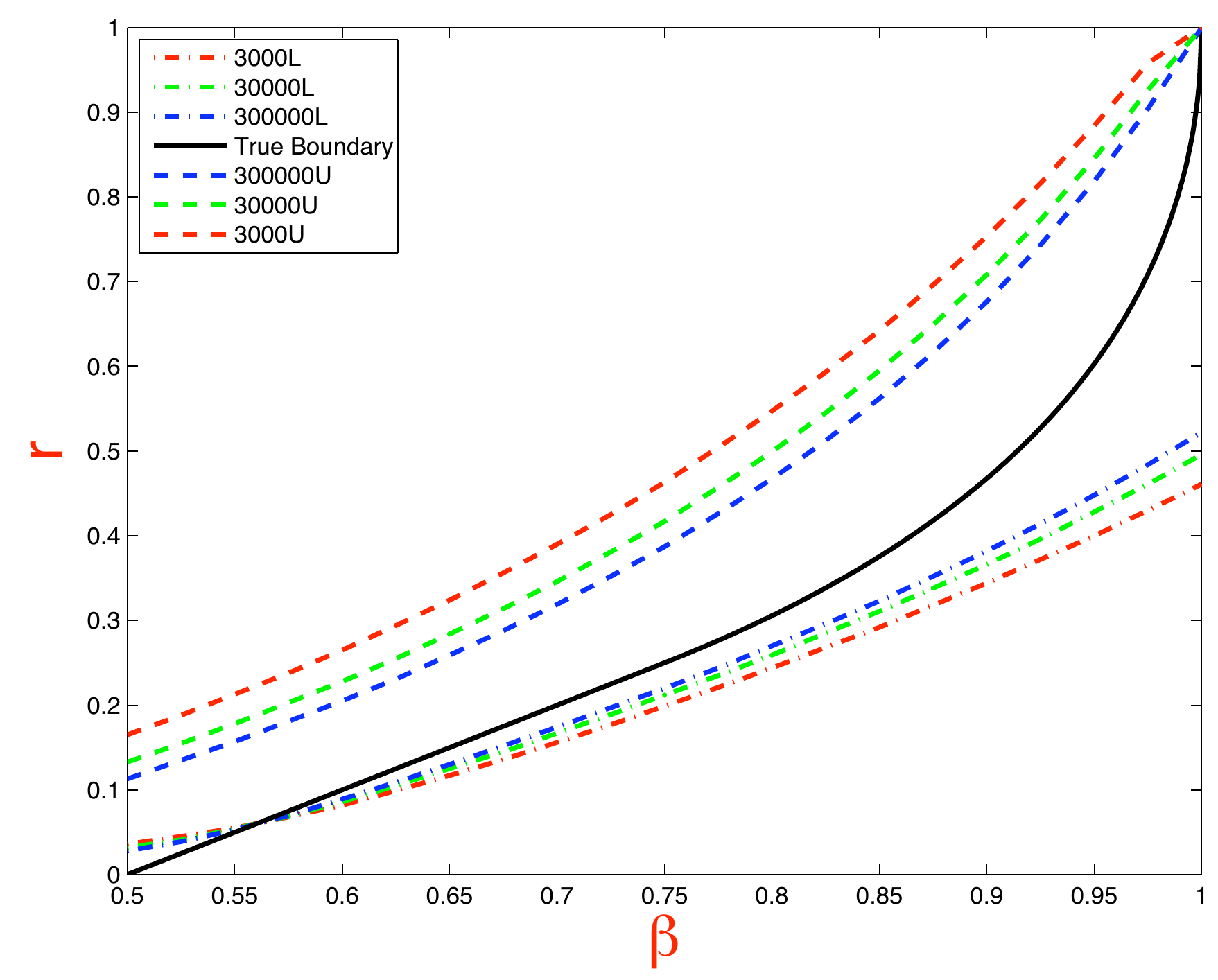}
\caption{Classification boundaries  for  $p = 3 \times 10^3 \times (1, 10, 100)$.  
For each $p$, both an upper boundary and a lower boundary are calculated which correspond to the 
proxy classification errors of $10\%$ and $40\%$ using ideal HCT. 
The gap between two boundaries gets smaller as $p$ increases. 
The solid black curve displays the asymptotic boundary $r = \rho^*(\beta)$.  }
\label{fig:finitep} 
\end{figure}

\subsection{Other work on HC}
\label{sec-OtherHC}
HC was originally proposed for use in a detection problem which 
has  nothing to do with threshold feature selection:
testing an intersection null hypothesis $\mu(j) = 0 \;  \forall j$ \cite{HC}.
The literature has developed since then.
Papers \cite{HJ1, HJ2} extends the optimality of  Higher Criticism  in detection to correlated setting.    
 Wellner  and his collaborators investigated the Higher
Criticism in the context of Goodness-of-fit,  see for example \cite{Wellner}.    Hall   and his 
collaborators   have investigated Higher Criticism for robustness, see for example \cite{Hall-robust}. HC has been
applied to data analysis in astronomy  (\cite{Cayon, Starck}) and
computational biology  \cite{Goeman}.

HC has been used as a principle to synthesize new procedures: Cai, Jin,
Low  in \cite{CJL}  (see also Meinshausen and Rice \cite{Rice}) use
Higher Criticism to motivate  estimators for the
proportion $\epsilon$ of non-null effects.  

\subsection{ HC in classification} 

Higher Criticism was previously applied to high-dimensional classification 
in Hall, Pittelkow, Ghosh (2008) \cite{Hall}, but there is  
a key conceptual difference from the present paper.      
Our approach  uses HC  in classifier design -- it selects features
in designing a linear classifier; but the actual classification
decisions are made by the linear classifier when presented with 
specific test feature vectors. The approach in \cite{Hall}  
uses HC to directly make decisions from feature vector data.  
Let's call these two different strategies {\it HC-based feature selection}
(used here) and {\it HC-based decision} (used in \cite{Hall}.)

In the ARW model, for a given classifier performance level,
the HC-based decision strategy  requires  much stronger feature contrasts
to reach that level than the HC-based feature selection strategy.
In this paper, we have shown that
HC-feature-selection requires useful features to
have contrasts exceeding
\[
  \sqrt{2 \rho^*(\beta)  \log p }/\sqrt{n} \cdot (1+o(1)), \qquad (n,p) \goto \infty.
\]
HC-decision requires  
useful features to have contrasts exceeding 
\[
  \sqrt{2 \rho^*(\beta) \log p } \cdot (1 + o(1)), \qquad p \goto \infty.
\]
Therefore, if the number of training samples $n > 1$,
HC-feature selection has an advantage.
Imagine that we have $n=36$ training samples; 
we will find that HC feature selection can asymptotically 
detect features roughly  
1/6 the strength of using HC directly for decision.

\end{document}